\documentclass[12pt]{article}
\usepackage{latexsym,amssymb,amsmath}
\begin{document}
\baselineskip=18pt

\newcommand{\la}{\langle}
\newcommand{\ra}{\rangle}
\newcommand{\psp}{\vspace{0.4cm}}
\newcommand{\pse}{\vspace{0.2cm}}
\newcommand{\ptl}{\partial}
\newcommand{\dlt}{\delta}
\newcommand{\sgm}{\sigma}
\newcommand{\al}{\alpha}
\newcommand{\be}{\beta}
\newcommand{\G}{\Gamma}
\newcommand{\gm}{\gamma}
\newcommand{\vs}{\varsigma}
\newcommand{\Lmd}{\Lambda}
\newcommand{\lmd}{\lambda}
\newcommand{\td}{\tilde}
\newcommand{\vf}{\varphi}
\newcommand{\for}{\mbox{for}}
\newcommand{\wt}{\mbox{wt}\:}
\newcommand{\der}{\mbox{Der}}
\newcommand{\ad}{\mbox{ad}\:}
\newcommand{\stl}{\stackrel}
\newcommand{\ol}{\overline}
\newcommand{\ul}{\underline}
\newcommand{\es}{\epsilon}
\newcommand{\dmd}{\diamond}
\newcommand{\clt}{\clubsuit}
\newcommand{\mbb}{\mathbb}
\newcommand{\llra}{\Longleftrightarrow}
\newcommand{\vt}{\vartheta}
\newcommand{\rta}{\rightarrow}
\newcommand{\ves}{\varepsilon}
\newcommand{\dg}{\dag}
\newcommand{\Dlt}{\Delta}

\begin{center}{\Large \bf Flag Partial Differential Equations and}\end{center}
\begin{center}{\Large \bf Representations of Lie Algebras}\end{center} \vspace{0.2cm}

\begin{center}{\large Xiaoping Xu}\footnote{Research Supported
 by China NSF 10431040}
\end{center}

\begin{center}{Institute of mathematics, Academy of Mathematics \&
System Sciences}\end{center}
\begin{center}{Chinese Academy of Sciences, Beijing 100080, P. R.
China}\end{center}

\vspace{2cm}

 \begin{center}{\Large\bf Abstract}\end{center}

\vspace{1cm} {\small In this paper, we solve the initial value
problems of variable-coefficient generalized wave equations
associated with trees and a large family of linear
constant-coefficient partial differential equation by algebraic
methods. Moreover, we find all the polynomial solutions for a
3-dimensional variable-coefficient flag partial differential
equation of any order, the linear wave equation with dissipation
and the generalized anisymmetrical Laplace equation. Furthermore,
the polynomial-trigonometric  solutions of  a generalized
Klein-Gordan equation associated with 3-dimensional generalized
Tricomi operator $\ptl_x^2+x\ptl_y^2+y\ptl_z^2$ are also given. As
applications to representations of Lie algebras, we find certain
irreducible polynomial  representations of the Lie algebras
$sl(n,\mbb{F}),\;so(n,\mbb{F})$ and the simple Lie algebra of type
$G_2$.}

 \vspace{0.8cm}

\section{Introduction}

Barros-Neto and Gel'fand [BG1,BG2] (1998, 2002) studied  solutions
of the equation
$$u_{xx}+xu_{yy}=\dlt(x-x_0,y-y_0)\eqno(1.1)$$
related to the Tricomi operator $\ptl_x^2+x\ptl_y^2$. A natural
generalization of the Tricomi operator is
$\ptl_{x_1}^2+x_1\ptl_{x_2}^2+\cdots+x_{n-1}\ptl_{x_n}^2$. The
equation
$$u_t=u_{x_1x_1}+u_{x_2x_2}+\cdots+u_{x_nx_n}\eqno(1.2)$$
is the well known classical heat conduction equation related to
the Laplacian operator
$\ptl_{x_1}^2+\ptl_{x_2}^2+\cdots+\ptl_{x_n}^2$. As pointed out in
[BG1, BG2], the Tricomi operator is an analogue of the Laplacian
operator. In [X2], we have solved the following natural analogue
of heat conduction equation:
$$u_t=u_{x_1x_1}+x_1u_{x_2x_2}+\cdots+x_{n-1}u_{x_nx_n}.\eqno(1.3)$$
Indeed we have solved analogous heat conduction equations related
to more general generalized Tricomi operators associated with tree
graphs.

A {\it tree} ${\cal T}$ consists of a finite set of {\it nodes}
${\cal N}= \{\iota_1,\iota_2,...,\iota_n\}$ and a set of {\it
edges}
$${\cal E}\subset\{(\iota_i,\iota_j)\mid 1\leq i<j\leq n\}\eqno(1.4)$$
such that for each node $\iota_i\in{\cal N}$, there exists a
unique sequence $\{\iota_{i_1},\iota_{i_2},...,\iota_{i_r}\}$ of
nodes with $1=i_1<i_2<\cdots<i_{r-1}<i_r=i$ for which
$$(\iota_{i_1},\iota_{i_2}),(\iota_{i_2},\iota_{i_3}),...,(\iota_{i_{r-2}},\iota_{i_{r-1}}),
(\iota_{i_{r-1}},\iota_{i_r})\in{\cal E}.\eqno(1.5)$$  We also
denote the tree ${\cal T}=({\cal N},{\cal E})$.
 For a tree ${\cal T}=({\cal N},{\cal E})$, we call the
differential operator
$$d_{\cal T}=\ptl_{x_1}^2+\sum_{(\iota_i,\iota_j)\in{\cal
E}}x_i\ptl_{x_j}^2\eqno(1.6)$$ a {\it generalized Tricomi operator
of type} ${\cal T}$. In [X2], we have solved the following
analogue of heat conduction equation:
$$u_t=d_{\cal T}(u)\eqno(1.7)$$
 subject to the initial condition:
$$u(0,x_1,...,x_n)=f(x_1,...,x_n)\qquad\for\;\;x_i\in[-a_i,a_i],\eqno(1.8)$$
where $f$ is a given continuous function and $a_i$ are given
positive real constants.

The equation
$$u_{tt}-u_{x_1x_1}-u_{x_2x_2}-\cdots-u_{x_nx_n}
=0\eqno(1.9)$$ is the well-known wave equation
 associated with the Laplace operator
$\ptl_{x_1}^2+\ptl_{x_2}^2+\cdots+\ptl_{x_n}^2$. One of our goals
in this paper is to solve the {\it generalized wave equation}
$$u_{tt}-d_{\cal T}(u)=0\eqno(1.10)$$
 subject to the initial conditions:
$$u(0,x_1,...,x_n)=g_0(x_1,...,x_n),\;\;\;u_t(0,x_1,...,x_n)=g_1(x_1,...,x_n)
\eqno(1.11)$$ for $x_i\in[-a_i,a_i]$, where $g_1,g_2$ are given
continuous functions and $a_i$ are given positive real constants.
This is done by using our results in [X2] and an algebraic method
of solving the following {\it differential equation of flag type}:
$$(d_1+f_1d_2+f_2d_3+\cdots+f_{n-1}d_n)(u)=0,\eqno(1.12)$$
where $d_1,d_2,...,d_n$ are certain commuting locally nilpotent
differential operators on the polynomial algebra
$\mbb{R}[x_1,x_2,...,x_n]$ and $f_1,...,f_{n-1}$ are polynomials
satisfying
$$d_i(f_j)=0\qquad\mbox{if}\;\;i>j.\eqno(1.13)$$
Another slightly adjusted method enables us to solve a large
family of linear constant-coefficient partial differential
equations subject to initial conditions easily. This family
contains almost all well-known linear constant-coefficient partial
differential equations such as the Helmholtz equation and the
Klein-Gordon equation. A general equation in this family can not
be solved by separation of variables.

Polynomial solutions of linear partial differential equations are
important for many reasons. For instance, polynomial solutions of
Laplace equations are called {\it harmonic polynomials}, which are
fundamental objects in differential geometry and analysis. In this
paper, we have used our algebraic methods to find a basis of the
space of all polynomial solutions for the following differential
equations: (1)
$$\ptl_x^{m_1}(u)+x^{n_1}\ptl_y^{m_2}(u)+y^{n_2}\ptl_z^{m_3}(u)=0,\eqno(1.14)$$
where $m_1,m_2,m_3$ are positive integers and $n_1,n_2$ are
nonnegative integers; (2) the linear wave equation with
dissipation:
$$u_{tt}+u_t-u_{x_1x_1}-u_{x_2x_2}-\cdots-u_{x_nx_n}=0;\eqno(1.15)$$
(3)
$$u_{tt}+\frac{\lmd}{t}u_t-\es(u_{x_1x_1}+u_{x_2x_2}+\cdots+u_{x_nx_n})=0,
 \eqno(1.16)$$ where $\es\in\{1,-1\}$
and  $\lmd$ is a nonzero real constant. When $m_1=m_2=m_3=2$ and
$n_1=n_2=0$, (1.14) is exactly the 3-dimensional Laplace equation.
If $\es=-1$, (1.16) is the generalized anisymmetrical Laplace
equation (cf. [A]). Moreover,  the special Euler-Poisson-Darboux
equation:
$$u_{tt}-u_{x_1x_1}-u_{x_2x_2}-\cdots-u_{x_nx_n}-\frac{m(m+1)}{t^2}u=0\eqno(1.17)$$
can be changed into a equation of type (1.16) with $\es=1$.
Furthermore, we have found the polynomial-trigonometric  solutions
of the generalized Klein-Gordan equation
$$u_{tt}-u_{xx}-xu_{yy}-yu_{zz}+a^2u=0\eqno(1.18)$$
associated with 3-dimensional generalized Tricomi operator
$\ptl_x^2+x\ptl_y^2+y\ptl_z^2$.

Flag partial differential equations also naturally appear in the
representation theory of Lie algebras. Let ${\cal G}$ be a
finite-dimensional simple Lie algebra with the Cartan root-space
decomposition:
$${\cal G}={\cal H}\oplus \bigoplus_{\al\in\Phi}{\cal
G}_\al.\eqno(1.19)$$ Suppose that $M$ is a finite-dimensional
${\cal G}$-module. Take a weight-vector basis $\{v_1,...,v_n\}$ of
$M$. Write
$$\xi(v_i)=\sum_{j=1}^n\rho_{i,j}(\xi)v_j\qquad\for\;\;\xi\in{\cal
G}.\eqno(1.20)$$ We define an action of ${\cal G}$ on
$\mbb{R}[x_1,...,x_n]$ by
$$\xi(g)=\sum_{i,j=1}^n\rho_{i,j}(\xi)x_j\ptl_{x_i}(g)\qquad\for\;\;\xi\in{\cal
G},\;g\in \mbb{R}[x_1,...,x_n].\eqno(1.21)$$ Then
$\mbb{R}[x_1,...,x_n]$ forms a ${\cal G}$-module isomorphic to the
symmetric tensor $S(M)$ over $M$. Take a base $\Pi$ of the root
system $\Phi$ and $0\neq \xi_\al\in {\cal G}_\al$ for $\al\in\Pi$.
According to Weyl's theorem of completely reducibility,
$\mbb{R}[x_1,...,x_n]$ can be decomposed as a direct sum of
irreducible ${\cal G}$-submodules. Such a decomposition is
completely determined by the polynomial solutions (singular
vectors or highest-weight vectors) of the following system of flag
partial differential equations:
$$\sum_{i,j=1}^n\rho_{i,j}(\xi_\al)x_j\ptl_{x_i}(u)=0,\qquad\al\in\Pi.
\eqno(1.22)$$

As simple applications of our algebraic methods of solving
differential equations, we find a basis for certain irreducible
polynomial representations of the Lie algebras
$sl(n,\mbb{F}),\;so(n,\mbb{F})$ and the simple Lie algebra of type
$G_2$.

In Section 2, we will present our algebraic methods  and find all
polynomial solutions of the above concerned partial differential
equations. Section 3 is devoted to solve two initial value
problems we mentioned in the above. In Section 4, certain
polynomial representations of Lie algebras will be given.

\section{Polynomial Solutions}

Take any subfield $\mbb{F}$ of the field $\mbb{C}$ of complex
numbers. We assume all the vector spaces are over $\mbb{F}$ unless
it is specified. For a positive integer $n$, we denote by
$${\cal A}=\mbb{F}[x_1,...,x_n],\eqno(2.1)$$
the algebra of polynomials in $n$ variables. The following
algebraic result is one of  our key lemmas of solving concerned
differential equations. \psp

{\bf Lemma 2.1}. {\it Suppose that ${\cal A}$ is free module over
a subalgebra ${\cal B}$ generated  by a filtrated subspace
$V=\bigcup_{r=0}^\infty V_r$ (i.e., $V_r\subset V_{r+1}$). Let
$T_1$ be a linear operator on ${\cal A}$ with a right inverse
$T_1^-$ such that
$$T_1({\cal B}),\;T_1^-({\cal B})\subset{\cal B},\qquad
T_1(\eta_1\eta_2)=T_1(\eta_1)\eta_2,\qquad
T_1^-(\eta_1\eta_2)=T_1^-(\eta_1)\eta_2 \eqno(2.2)$$ for $\eta_1
\in {\cal B},\;\eta_2\in V$, and let $T_2$ be a linear operator on
${\cal A}$ such that
$$ T_2(V_{r+1})\subset {\cal B}V_r,\;\;
T_2(f\zeta)=fT_2(\zeta) \qquad\for\;\;0\leq
r\in\mbb{Z},\;\;f\in{\cal B},\;\zeta\in{\cal A}.\eqno(2.3)$$ Then
we
have \begin{eqnarray*}\hspace{1cm}&&\{g\in{\cal A}\mid (T_1+T_2)(g)=0\}\\
& =&\mbox{Span}\{ \sum_{i=0}^\infty(-T_1^-T_2)^i(hg)\mid g\in
V,\;h\in {\cal B};\;T_1(h)=0\}, \hspace{3.6cm}(2.4)\end{eqnarray*}
where the summation is finite under our assumption. Moreover, the
operator $\sum_{i=0}^\infty(-T_1^-T_2)^iT_1^-$ is a right inverse
of $T_1+T_2$}.

{\it Proof}. For $h\in {\cal B}$ such that $T_1(h)=0$ and $g\in
V$, we have
\begin{eqnarray*}& &(T_1+T_2)(\sum_{i=0}^\infty(-T_1^-T_2)^i(hg))\\
&=&T_1(hg)-\sum_{i=1}^\infty T_1[T_1^-T_2(-T_1^-T_2)^{i-1}(hg)]+
\sum_{i=0}^\infty T_2[(-T_1^-)^i(hg)]\\
&=&T_1(h)g-\sum_{i=1}^\infty (T_1T_1^-)T_2(-T_1^-T_2)^{i-1}(hg)+
\sum_{i=0}^\infty T_2(-T_1^-T_2)^i(hg)\\
&=&-\sum_{i=1}^\infty T_2(-T_1^-T_2)^{i-1}(hg)+ \sum_{i=0}^\infty
T_2(-T_1^-T_2)^i(hg)=0\hspace{4.8cm} (2.5)\end{eqnarray*} by
(2.2).  Denote by $\mbb{N}$ the set of nonnegative integers. Set
$V_{-1}=\{0\}$. For $k\in\mbb{N}$, we take $\{\psi_i\mid i\in
I_k\}\subset V_k$ such that
$$\{\psi_i+V_{k-1}\mid i\in I_k\}\;\;\mbox{forms a basis
of}\;\;V_k/V_{k-1},\eqno(2.6)$$ where $I_k$ is an index set.
 Let
$${\cal A}^{(m)}={\cal B}V_m=\sum_{s=0}^m\;\sum_{i\in I_s}{\cal
B}\psi_i.\eqno(2.7)$$ Obviously,
$$T_1({\cal
A}^{(m)}),\;T_1^-({\cal A}^{(m)}),\; T_2({\cal A}^{(m+1)})\subset
{\cal A}^{(m)}\qquad\for\;\;m\in\mbb{N}\eqno(2.8)$$ by (2.2) and
(2.3), and
$${\cal A}=\bigcup_{m=0}^\infty {\cal A}^{(m)}.\eqno(2.9)$$

  Suppose $\phi\in {\cal
A}^{(m)}$ such that $(T_1+T_2)(\phi)=0$. If $m=0$, then
$$\phi=\sum_{i\in I_0}h_i\psi_i,\qquad h_i\in {\cal
B}.\eqno(2.10)$$ Now
$$0=(T_1+T_2)(\phi)=\sum_{i\in I_0}T_1(h_i)\psi_i+\sum_{i\in
I_0}h_iT_2(\psi_i)=\sum_{i\in I_0}T_1(h_i)\psi_i,\eqno(2.11)$$
Since $T_1(h_i)\in {\cal B}$ by (2.2) and ${\cal A}$ is a free
${\cal B}$-module generated by $V$, we have $T_1(h_i)=0$ for $i\in
I_0$. Denote by ${\cal S}$ the right hand side of the equation
(2.4). Then
$$\phi=\sum_{i\in I_0}\;\sum_{m=0}^\infty(-T_1^-T_2)^m(h_i\psi_i)\in {\cal
S}.\eqno(2.12)$$

Suppose $m>0$. We write
$$\phi=\sum_{i\in I_m}h_i\psi_i+\phi',\qquad h_i\in{\cal
B},\;\phi'\in {\cal A}^{(m-1)}.\eqno(2.13)$$ Then
$$0=(T_1+T_2)(\phi)=\sum_{i\in
I_m}T_1(h_i)\psi_i+T_1(\phi')+T_2(\phi).\eqno(2.14)$$ Since
$T_1(\phi')+T_2(\phi)\in {\cal A}^{(m-1)}$, we have $T_1(h_i)=0$
for $i\in I_m$. Now
$$\phi-\sum_{i\in
I_m}\sum_{k=0}^\infty(-T_1^-T_2)^k(h\psi_i)=\phi'-\sum_{i\in
I_m}\sum_{k=1}^\infty(-T_1^-T_2)^k(h_i\psi_i)\in {\cal
A}^{(m-1)}\eqno(2.15)$$ and (2.5) implies
$$(T_1+T_2)(\phi-\sum_{i\in
I_m}\sum_{k=0}^\infty(-T_1^-T_2)^k(h_i\psi_i))=0.\eqno(2.16)$$ By
induction on $m$, $$\phi-\sum_{i\in
I_m}\sum_{k=0}^\infty(-T_1^-T_2)^k(h_i\psi_i)\in {\cal
S}.\eqno(2.17)$$ Therefore, $\phi\in {\cal S}.$

For any $f\in {\cal A}$, we have:
\begin{eqnarray*} \hspace{1cm} & &(T_1+T_2)(\sum_{i=0}^\infty(-T_1^-T_2)^iT_1^-)(f)
\\ &=&f-\sum_{i=1}^\infty
T_2(-T_1^-T_2)^{i-1}T_1^-(f)+\sum_{i=0}^\infty
T_2(-T_1^-T_2)^iT_1^-(f)=f.\hspace{2.5cm}(2.18)\end{eqnarray*}
Thus
 the operator
$\sum_{i=0}^\infty(-T_1^-T_2)^iT_1^-$ is a right inverse of
$T_1+T_2.\qquad\Box$\psp

 We remark that the above operator $T_1$ and $T_2$ may not
 commute.

For convenience, we take the following notation of indices:
$$\ol{i,j}=\{i,i+1,...,j\},\eqno(2.19)$$
where $i\leq j$ are integers. Define
$$x^\al=x_1^{\al_1}x_2^{\al_2}\cdots
x_n^{\al_n}\qquad\for\;\;\al=(\al_1,...,\al_n)\in\mbb{N}^{\:n}.\eqno(2.20)$$
Moreover, we denote
$$\es_i=(0,...,0,\stl{i}{1},0,...,0)\in \mbb{N}^{\:n}.\eqno(2.21)$$
 For each
$i\in\ol{1,n}$, we define the linear operator $\int_{(x_i)}$ on
${\cal A}$ by:
$$\int_{(x_i)}(x^\al)=\frac{x^{\al+\es_i}}{\al_i+1}\qquad\for\;\;\al\in
\mbb{N}^{\:n}.\eqno(2.22)$$ Furthermore, we let
$$\int_{(x_i)}^{(0)}=1,\qquad\int_{(x_i)}^{(m)}=\stl{m}{\overbrace{\int_{(x_i)}\cdots\int_{(x_i)}}}
\qquad\for\; \;0<m\in\mbb{Z}\eqno(2.23)$$ and denote
$$\ptl^{\al}=\ptl_{x_1}^{\al_1}\ptl_{x_2}^{\al_2}\cdots
\ptl_{x_n}^{\al_n},\;\;
\int^{(\al)}=\int_{(x_1)}^{(\al_1)}\int_{(x_2)}^{(\al_2)}\cdots
\int_{(x_n)}^{(\al_n)}\qquad\for\;\;\al\in
\mbb{N}^{\:n}.\eqno(2.24)$$ Obviously, $\int^{(\al)}$ is a right
inverse of $\ptl^\al$ for $\al\in \mbb{N}^{\:n}.$ We remark that
$\int^{(\al)}\ptl^\al\neq 1$ if $\al\neq 0$ due to
$\ptl^\al(1)=0$.

Consider the wave equation in Riemannian space with a nontrivial
conformal group (cf. [I1]):
$$u_{tt}-u_{x_1x_1}-\sum_{i,j=2}^ng_{i,j}(x_1-t)u_{x_ix_j}=0,\eqno(2.25)$$
where we assume that $g_{i,j}(z)$ are one-variable polynomials.
Change variables:
$$z_0=x_1+t,\qquad z_1=x_1-t.\eqno(2.26)$$
Then
$$\ptl_t^2=(\ptl_{z_0}-\ptl_{z_1})^2,\qquad
\ptl_{x_1}^2=(\ptl_{z_0}+\ptl_{z_1})^2.\eqno(2.27)$$ So the
equation (2.25) changes to:
$$2\ptl_{z_0}\ptl_{z_1}+
\sum_{i,j=2}^ng_{i,j}(z_1)u_{x_ix_j}=0.\eqno(2.28)$$ Denote
$$T_1=2\ptl_{z_0}\ptl_{z_1},\qquad
T_2=\sum_{i,j=2}^ng_{i,j}(z_1)\ptl_{x_i}\ptl_{x_j}.\eqno(2.29)$$
Take $T_1^-=\frac{1}{2}\int_{(z_0)}\int_{(z_1)}$, and
$${\cal B}=\mbb{F}[z_0,z_1],\qquad V=\mbb{F}[x_2,...,x_n],\qquad
V_r=\{f\in V\mid\mbox{deg}\;f\leq r\}.\eqno(2.30)$$ Then the
conditions in Lemma 2.1 hold. Thus we have:
 \psp

{\bf Theorem 2.2}. {\it The space of all polynomial solutions for
the equation (2.24) is:
\begin{eqnarray*} \hspace{2cm}& &\mbox{Span}\:\{\sum_{m=0}^\infty(-2)^{-m}
(\sum_{i,j=2}^n\int_{(z_0)}\int_{(z_1)}g_{i,j}(z_1)\ptl_{x_i}\ptl_{x_j})^m(f_0g_0+f_1g_1)\\
& & \mid
f_0\in\mbb{F}[z_0],\;f_1\in\mbb{F}[z_1],\;g_0,g_1\in\mbb{F}[x_2,...,x_n]\}
\hspace{4.3cm}(2.31)\end{eqnarray*} with $z_0,z_1$ defined in
(2.26).}\psp

Let $m_1,m_2,...,m_n$ be positive integers. According to Lemma
2.1, the set
\begin{eqnarray*}\hspace{1.9cm}& &\{\sum_{k_2,...,k_n=0}^\infty(-1)^{k_2+\cdots+k_n}{k_2+\cdots+k_k\choose
k_2,...,k_n} \int_{(x_1)}^{((k_2+\cdots +k_n)m_1)}(x_1^{\ell_1})\\
& &\times\ptl_{x_2}^{k_2m_2}(x_2^{\ell_2})\cdots
\ptl_{x_n}^{k_nm_n}(x_n^{\ell_n})\mid
\ell_1\in\ol{0,m_1-1},\;\ell_2,...,\ell_n\in\mbb{N}\}\hspace{2cm}
(2.32)\end{eqnarray*} forms a basis of the space of polynomial
solutions for the equation
$$(\ptl_{x_1}^{m_1}+\ptl_{x_2}^{m_2}+\cdots+\ptl_{x_n}^{m_n})(u)=0\eqno(2.33)$$
in  ${\cal A}$. In particular,
\begin{eqnarray*}\hspace{2cm}& &\{\sum_{r_2,...,r_n=0}^{\infty}\frac{(-1)^{r_2+\cdots+
r_k}{r_2+\cdots+r_n\choose r_2,...,r_n}\prod_{i=2}^n{\ell_i\choose
2r_i}}{(1+2\es(r_2+\cdots+r_n)){2(r_2+\cdots+r_n)\choose
2r_2,...,2r_n}} x_1^{\es+2(r_2+\cdots +r_n)}\\ & &\times
x_2^{\ell_2-2r_2}\cdots x_n^{\ell_n-2r_n}
\mid\es\in\{0,1\};\;\ell_2,...,\ell_n\in\mbb{N}\}\hspace{4.1cm}(2.34)\end{eqnarray*}
forms a basis of the space of harmonic polynomials in $n$
variables, that is, the space of polynomial solutions for the
Laplace equation:
$$u_{x_1x_1}+u_{x_2x_2}+\cdots+u_{x_nx_n}=0.\eqno(2.35)$$

The above results can theoretically generalized as follows. Let
$$f_i\in\mbb{F}[x_1,...,x_i]\qquad \for\;\;i\in\ol{1,n-1}.\eqno(2.36)$$
Consider the equation:
$$(\ptl_{x_1}^{m_1}+f_1\ptl_{x_2}^{m_2}+\cdots+f_{n-1}\ptl_{x_n}^{m_n})(u)=0
\eqno(2.37)$$
 Denote
$$d_1=\ptl_{x_1}^{m_1},\;\;
d_r=\ptl_{x_1}^{m_1}+f_1\ptl_{x_2}^{m_2}+\cdots+f_{r-1}\ptl_{x_r}^{m_r}\qquad\for\;\;r
\in\ol{2,n}.\eqno(2.38)$$ We will apply Lemma 2.1 with
$T_1=d_r,\;T_2=\sum_{i=r}^{n-1}f_i\ptl_{x_{i+1}}^{m_{i+1}}$ and
${\cal B}=\mbb{F}[x_1,...,x_r],\;V=\mbb{F}[x_{r+1},...,x_n]$,
$$V_k=\mbox{Span}\:\{x_{r+1}^{\ell_{r+1}}\cdots x_n^{\ell_n}\mid
\ell_s\in\mbb{N},\;\ell_{r+1}+\sum_{i=r+2}^n\ell_i(\mbox{deg}\:f_{r+1}+1)\cdots
(\mbox{deg}\:f_{i-1}+1)\leq k\}.\eqno(2.39)$$
 Take a right inverse
$d_1^-=\int_{(x_1)}^{(m_1)}$.  Suppose that we have found a right
inverse $d_s^-$ of $d_s$ for some $s\in\ol{1,n-1}$ such that
$$x_id_s^-=d_s^-x_i,\;\;
\ptl_{x_i}d_s^-=d_s^-\ptl_{x_i}\qquad\for\;\;i\in\ol{s+1,n}.\eqno(2.40)$$
Lemma 2.1 enable us to take
$$d_{s+1}^-=\sum_{i=0}^\infty(-d_s^-f_s)^id_s^-\ptl_{x_{s+1}}^{im_{s+1}}\eqno(2.41)$$
as a right inverse of $d_{s+1}$.  Obviously, $$
x_id_{s+1}^-=d_{s+1}^-x_i,\;\;
\ptl_{x_i}d_{s+1}^-=d_{s+1}^-\ptl_{x_i}\qquad\for\;\;i\in\ol{s+2,n}
\eqno(2.42)$$ according to (2.38). By induction, we have found a
right inverse $d_s^-$ of $d_s$ such that (2.40) holds for each
$s\in\ol{1,n}$.

We set
$${\cal S}_r=\{g\in \mbb{F}[x_1,...,x_r]\mid
d_r(g)=0\}\qquad\for\;\;r\in\ol{1,k}.\eqno(2.43)$$ By (2.38),
$${\cal S}_1=\sum_{i=0}^{m_1-1}\mbb{F}x_1^i.\eqno(2.44)$$
Suppose that we have found ${\cal S}_r$ for some $r\in
\ol{1,n-1}$. Given $h\in {\cal S}_r$ and $\ell\in \mbb{N}$, we
define
$$\sgm_{r+1,\ell}(h)=\sum_{i=0}^\infty
(-d_r^-f_r)^i(h)\ptl_{x_{r+1}}^{im_{r+1}}(x_{r+1}^{\ell}),\eqno(2.45)$$
which is actually a finite summation. Lemma 2.1 says
$${\cal S}_{r+1}=\sum_{\ell=0}^\infty \sgm_{r+1,\ell}({\cal
S}_r).\eqno(2.46)$$ By  induction, we obtain:\psp

{\bf Theorem 2.3}. {\it The set
$$\{\sgm_{n,\ell_n}\sgm_{n-1,\ell_{n-1}}\cdots\sgm_{2,\ell_2}(x_1^{\ell_1})\mid
\ell_1\in\ol{0,m_1-1},\;\ell_2,...,\ell_n\in\mbb{N}\}\eqno(2.47)$$
forms a basis of the polynomial solution space ${\cal S}_n$ of the
partial differential equation (2.37).}\psp

{\bf Example 2.1}. Let $m_1,m_2,m_3,n_1,n_2$ be positive integers.
Consider the following equations
$$\ptl_{x}^{m_1}(u)+x^{n_1}\ptl_{y}^{m_2}(u)+y^{n_2}\ptl_{z}^{m_3}(u)=0\eqno(2.48)$$
Now
$$d_1=\ptl_{x}^{m_1},\;\;\;d_2=\ptl_x^{m_1}+x^{n_1}\ptl_y^{n_2}.\eqno(2.49)$$
Take $d_1^-=\int_{(x)}^{(m_1)}$. Then
$$\sgm_{2,\ell_2}(x^{\ell_1})=\sum_{i=0}^\infty(-\int_{(x)}^{(m_1)}x^{n_1})^i
(x^{\ell_1})\ptl_{y}^{im_2}(y^{\ell_2}). \eqno(2.50)$$ Moreover,
$$d_2^-=\sum_{i=0}^\infty(-\int_{(x)}^{(m_1)}x^{n_1})^i\int_{(x)}^{(m_1)}\ptl_{y}^{im_2}
\eqno(2.51)$$ by Lemma 2.1. Hence
\begin{eqnarray*}(d_2^-y^{n_2})^k&=&\sum_{i_1,...,i_k=0}^\infty(-1)^{i_1+\cdots+i_k}
(\int_{(x)}^{(m_1)}x^{n_1})^{i_1}\int_{(x)}^{m_1}
(\int_{x}^{(m_1)}x^{n_1})^{i_2}\int_{(x)}^{(m_1)}\cdots
\\ &
&(\int_{(x)}^{(m_1)}x^{n_1})^{i_k}\int_{(x)}^{(m_1)}
\ptl_{y}^{i_1m_2}y^{n_2}\ptl_{y}^{i_2m_2}y^{n_2}\cdots\ptl_{y}^{i_km_2}y^{n_2}.
\hspace{3.7cm}(2.52)\end{eqnarray*} Thus
\begin{eqnarray*}& &\sgm_{3,\ell_3}\sgm_{2,\ell_2}(x^{\ell_1})\\&=&
\sum_{i=0}^\infty(-\int_{(x)}^{(m_1)}x^{n_1})^i
(x^{\ell_1})\ptl_{y}^{im_2}(y^{\ell_2})z^{\ell_3}+
\sum_{k=1}^\infty\sum_{i,i_1,...,i_k=0}^\infty(-1)^{k+i+i_1+\cdots+i_k}
\\ & &[(\int_{(x)}^{(m_1)}x^{n_1})^{i_1}\int_{(x)}^{(m_1)}
(\int_{(x)}^{(m_1)}x^{n_1})^{i_2}\int_{(x)}^{(m_1)}\cdots
(\int_{(x)}^{(m_1)}x^{n_1})^{i_k}\int_{(x)}^{(m_1)}(\int_{(x)}^{(m_1)}x^{n_1})^i]
(x^{\ell_1})\\ & &\times
[\ptl_{y}^{i_1m_2}y^{n_2}\ptl_{y}^{i_2m_2}y^{n_2}\cdots\ptl_{y}^{i_km_2}
y^{n_2}\ptl_{y}^{im_2}](y^{\ell_2})\ptl_{z}^{km_3}(z^{\ell_3})\\
 &=&\sum_{i=0}^\infty(-1)^i\frac{\prod_{s=0}^{im_2-1}(\ell_2-s)}{\prod_{r=0}^{i-1}\prod_{s=1}^{m_1}
 (s+n_1+\ell_1+r(m_1+n_1))}x^{\ell_1+i(m_1+n_1)}y^{\ell_2-im_2}\\
 &
 &+
 \sum_{k=1}^\infty\sum_{i_0,i_1,...,i_k=0}^\infty\frac{(-1)^{k+i+i_1+\cdots+i_k}}{\prod_{r_0=0}^{i_0-1}\prod_{s_0=1}^{m_1}
 (s+n_1+\ell_1+r_0(m_1+n_1))}\\ & &\times\frac{\left(\prod_{p=0}^k\sum_{r=0}^{i_pm_2-1}(pn_2+\ell_2-(\sum_{s=0}^{p-1}i_s)m_2-r)\right)
 \left(\prod_{r=0}^{km_3-1}(\ell_3-r)\right)}{\prod_{p=1}^k
\left(\prod_{r_p=0}^{i_p}\prod_{s=1}^{m_1}
 (s+\ell_1+(p-1)m_1+(r_p+\sum_{q=0}^{p-1}i_q)(m_1+n_1))\right)}\\
 & &\times
 x^{(\sum_{r=0}^ki_r)(m_1+n_1)+km_1+\ell_1}y^{kn_2+\ell_2-(\sum_{r=0}^ki_r)m_2}z^{\ell_3-km_3}
. \hspace{4.9cm}(2.53)\end{eqnarray*} \pse

In order to solve the linear wave equation with dissipation and
the generalized Klein-Gordan equation, we need the following
lemma. \psp

{\bf Lemma 2.4}. {\it Let $d=a\ptl_t+\ptl_t^2$ with $0\neq
a\in\mbb{F}$. Take a right inverse
$$d^-=\int_{(t)}\sum_{r=0}^\infty a^{-r-1}(-\ptl_t)^r\eqno(2.54)$$
of $d$. Then}
$$(d^-)^i(1)=\frac{t^i}{i!a^i}-\frac{t^{i-1}}{(i-2)!a^{i+1}}+\sum_{r=2}^{i-1}\frac{(-1)^r\prod_{s=1}^{r-1}(i+s)}
{(i-r-1)!r!a^{r+i}}t^{i-r}.\eqno(2.55)$$

{\it Proof}. For
$$f(t)=\sum_{i=1}^mb_it^i\in\mbb{F}[t]t,\eqno(2.56)$$
we have
$$d(f(t))=amb_mt^{m-1}+\sum_{i=1}^{m-1}i(ab_i+(i+1)b_{i+1})t^{i-1}.\eqno(2.57)$$
Thus $d(f(t))=0$ if and only if $f(t)\equiv 0$. So for any given
$g(t)\in\mbb{F}[t]$, there exists a unique $f(t)\in \mbb{F}[t]t$
such that $d(f(t))=g(t)$.

Set
$$\xi_{a,i}(t)=\frac{t^i}{i!a^i}-\frac{t^{i-1}}{(i-2)!a^{i+1}}+\sum_{r=2}^{i-1}\frac{(-1)^r\prod_{s=1}^{r-1}(i+s)}
{(i-r-1)!r!a^{r+i}}t^{i-r},\eqno(2.58)$$ where we treat
$$\xi_{a,0}(t)=1,\;\;\xi_{a,1}(t)=\frac{t}{a},\;\;\xi_{a,2}(t)=\frac{t^2}{2a^2}-\frac{t}{a^3}.
\eqno(2.59)$$ Easily verify $d(\xi_{a,i}(t))=\xi_{a,i-1}(t)$ for
$i=1,2$.

 Assume $i>2$. We have
\begin{eqnarray*} &&d(\xi_{a,i}(t))\\ &=&(a\ptl_t+\ptl_t^2)\left(\frac{t^i}{i!a^i}-
\frac{t^{i-1}}{(i-2)!a^{i+1}}+\sum_{r=2}^{i-1}\frac{(-1)^r\prod_{s=1}^{r-1}(i+s)}
{(i-r-1)!r!a^{r+i}}t^{i-r}\right)\\
&=&\frac{t^{i-1}}{(i-1)!a^{i-1}}-\frac{(i-1)t^{i-2}}{(i-2)!a^i}+
\sum_{r=2}^{i-1}\frac{(-1)^r(i-r)\prod_{s=1}^{r-1}(i+s)}
{(i-r-1)!r!a^{r+i-1}}t^{i-r-1}\\ &
&+\frac{t^{i-2}}{(i-2)!a^i}-\frac{(i-1)t^{i-3}}{(i-3)!a^{i+1}}+\sum_{r=2}^{i-1}
\frac{(-1)^r(i-r)\prod_{s=1}^{r-1}(i+s)}
{(i-r-2)!r!a^{r+i}}t^{i-r-2}\\
&=&\frac{t^{i-1}}{(i-1)!a^{i-1}}-\frac{t^{i-2}}{(i-3)!a^i}+
\frac{(i-2)(i+1)}
{(i-3)!2!a^{i+1}}t^{i-3}-\frac{(i-1)t^{i-3}}{(i-3)!a^{i+1}}\\ & &+
\sum_{r=3}^{i-1}(-1)^r\left[\frac{(i-r)\prod_{s=1}^{r-1}(i+s)}
{r!}-\frac{(i-r+1)\prod_{s=1}^{r-2}(i+s)}
{(r-1)!}\right]\frac{t^{i-r-1}}{(i-r-1)!a^{r+i-1}}\\ &=&
\frac{t^{i-1}}{(i-1)!a^{i-1}}-\frac{t^{i-2}}{(i-3)!a^i}+
\frac{(i-2)(i+1)-2(i-1)}{(i-3)!2!a^{i+1}}t^{i-3}
\\ & &+\sum_{r=3}^{i-1}(-1)^r\frac{[(i-r)(i+r-1)-r(i-r+1)]\prod_{s=1}^{r-2}(i+s)}
{(i-r-1)!r!a^{r+i-1}}t^{i-r-1}
\\ &=&
\frac{t^{i-1}}{(i-1)!a^{i-1}}-\frac{t^{i-2}}{(i-3)!a^i}+
\frac{i(i-3)}{(i-3)!2!a^{i+1}}t^{i-3}
\\ & &+\sum_{r=3}^{i-1}(-1)^r\frac{i(i-1-r)\prod_{s=1}^{r-2}(i+s)}
{(i-r-1)!r!a^{r+i-1}}t^{i-r-1}\\&=&
\frac{t^{i-1}}{(i-1)!a^{i-1}}-\frac{t^{i-2}}{(i-3)!a^i}+
\frac{i}{(i-4)!2!a^{i+1}}t^{i-3}
+\sum_{r=3}^{i-2}(-1)^r\frac{i\prod_{s=1}^{r-2}(i+s)}
{(i-r-2)!r!a^{r+i-1}}t^{i-r-1}
\\&=&\frac{t^{i-1}}{(i-1)!a^{i-1}}-\frac{t^{i-2}}{(i-3)!a^i}+
\sum_{r=2}^{i-2}(-1)^r\frac{\prod_{s=1}^{r-1}(i-1+s)}
{(i-r-2)!r!a^{r+i-1}}t^{i-r-1}\\ &
=&\xi_{a,i-1}(t).\hspace{12.5cm}(2.60)\end{eqnarray*} Since
$(d^-)^0(1)=1,\;(d^-)^i(1)\in\mbb{F}[t]t$ by (2.54) and
$d[(d^-)^i(1)]=(d^-)^{i-1}(1)$ for $i\in\mbb{N}+1$, we have
$(d^-)^r(1)=\xi_{a,r}(t)$ for $r\in\mbb{N}$ by the uniqueness,
that is, (2.55) holds.$\qquad\Box$\psp

By Lemma 2.1 and the above lemma, we obtain:\psp

{\bf Theorem 2.5}. {\it The set
\begin{eqnarray*} \hspace{2cm}& &\{\sum_{r_1,...,r_n=0}^{\infty}{r_1+\cdots+r_n\choose
r_1,...,r_n}\left[\prod_{i=1}^n(2r_i)!{\ell_i\choose 2r_i}\right]
\\ & &\times\xi_{1,r_1+\cdots +r_n}(t)x_1^{\ell_1-2r_1}\cdots
x_n^{\ell_n-2r_n} \mid
\ell_1,...,\ell_n\in\mbb{N}\}\hspace{3.8cm}(2.61)\end{eqnarray*}forms
a basis of the space of polynomial solution for the linear wave
equation with dissipation:}
$$u_{tt}+u_t-u_{x_1x_1}-u_{x_2x_2}-\cdots-u_{x_nx_n}=0\eqno(2.62)$$
\pse

 Consider the following generalization of
Klein-Gordan equation:
$$u_{tt}-u_{xx}-xu_{yy}-yu_{zz}+a^2u=0,\eqno(2.63)$$ where $a$ is a nonzero
real number. Changing variable $u=e^{a\sqrt{-1} t}v$, we get
$$v_{tt}+2a\sqrt{-1}v_t-v_{xx}-xv_{yy}-yv_{zz}=0.\eqno(2.64)$$ We write
$$\xi_{2a
\sqrt{-1},i}=\zeta_{i,0}(t)+\zeta_{i,1}(t)\sqrt{-1},\eqno(2.65)$$
where $\zeta_{i,0}(t)$ and $\zeta_{i,1}(t)$ are real functions.
According to (2.58),
$$\zeta_{2i,0}(t)=(-1)^i\left[\frac{t^{2i}}{(2i)!(2a)^{2i}}
+\sum_{r=1}^{i-1}\frac{(-1)^r\prod_{s=1}^{2r-1}(2i+s)}{(2r)!(2(i-r)-1)!(2
a)^{2(i+r)}}t^{2(i-r)}\right],\eqno(2.66)$$
\begin{eqnarray*}\hspace{2cm}& &\zeta_{2r,1}(t)=(-1)^i[\frac{t^{2r-1}}{(2r-2)!
(2a)^{2i+1}}\\ & &
+\sum_{r=1}^{i-1}\frac{(-1)^r\prod_{s=1}^{2r}(2i+s)}{(2r+1)!
[2(i-r-1)]!(2a)^{2i+2r+1}}t^{2i-2r-1}],\hspace{3.5cm}(2.67)\end{eqnarray*}
\begin{eqnarray*}\hspace{2cm}& &
\zeta_{2i+1,0}(t)=(-1)^i[\frac{t^{2i}}{(2i-1)! (2a)^{2(i+1)}}\\ &
&
+\sum_{r=1}^{i-1}\frac{(-1)^r\prod_{s=1}^{2r}(2i+s+1)}{(2r+1)!(2i-2r-1)!
(2a)^{2(i+r+1)}}t^{2(i-r)}],\hspace{3.7cm}(2.68)\end{eqnarray*}
$$
\zeta_{2r+1,1}(t)=(-1)^{i+1}
\left[\frac{t^{2i+1}}{(2i+1)!(2a)^{2i+1}}+
\sum_{r=1}^i\frac{(-1)^r\prod_{s=1}^{2r-1}(2i+s+1)}{(2r)!(2i-2r)!(2
a)^{2i+2r+1}}t^{2i-2r+1}\right].\eqno(2.69)$$

On the other hand, (5.22) in [X2] tells
$$e^{t(\ptl_{x}^2+x\ptl_{y}^2+y\ptl_{z}^2)}=e^{\xi_3}e^{\xi_2}
e^{\xi_1}\eqno(2.70)$$ with
\begin{eqnarray*}\hspace{1.6cm} \xi_1&=&
\int_0^t(\ptl_{x}+\int_0^{y_1}(\ptl_{y}
+y_2\ptl_{z}^2)^2dy_2)^2dy_1\\ &=&
\int_0^t(\ptl_{x}+y_1\ptl_{y}^2+y_1^2\ptl_{y}\ptl_{z}^2+y_1^3\ptl_{z}^4/3)^2dy_1\\
&=&t\ptl_{x}^2+t^2\ptl_{x}\ptl_{y}^2+\frac{t^3}{3}(\ptl_{y}^4+2\ptl_{x}\ptl_{y}
\ptl_{z}^2)+\frac{t^4}{6}(3\ptl_{y}^3\ptl_{z}^2+\ptl_{x}\ptl_{z}^4)\\
& & +\frac{t^5}{3} \ptl_{y}^2\ptl_{z}^4+\frac{t^6}{9}
\ptl_{y}\ptl_{z}^6+\frac{t^7}{21}
\ptl_{z}^8,\hspace{7.1cm}(2.71)\end{eqnarray*}
$$\xi_2=x(t\ptl_{y}^2+t^2\ptl_{y}\ptl_{z}^2+t^3\ptl_{z}^4/3),\;\;\;
\xi_3=ty\ptl_{z}^2.\eqno(2.72)$$ For
$\al=(\al_1,...,\al_{11})\in\mbb{N}^{11}$, we define
$$|\al|_1=\sum_{i=1}^7i\al_i+\al_8+2\al_9+3\al_{10}+\al_{11},\;\;\al!=\prod_{i=1}^{11}\al_i!
\eqno(2.73)$$ and
\begin{eqnarray*}\hspace{1cm}d_\al&=&
\frac{|\al|_1!}{\al!}2^{-\al_4}3^{-\al_3-\al_4-\al_5-2\al_6-\al_7-\al_{10}}7^{-\al_7}
x^{\al_8+\al_9+\al_{10}} y^{\al_{11}}\\
&&\times\ptl_{x}^{2\al_1+\al_2}
\ptl_{y}^{2(\al_2+\al_5+\al_8)+\al_6+\al_9}
\ptl_{z}^{2(2\al_5+3\al_6+4\al_7+\al_9+2\al_{10}+ \al_{11})}\\ &
&\times(\ptl_{y}^4+2\ptl_{x}\ptl_{y}
\ptl_{z}^2)^{\al_3}(3\ptl_{y}^3\ptl_{z}^2+\ptl_{x}\ptl_{z}^4)^{\al_4}.
\hspace{6cm}(2.74)\end{eqnarray*} Then
$$\sum_{m_1,m_2,m_3=0}^\infty\mbb{C}
(\sum_{\al\in\mbb{N}^{11}}\zeta_{|\al|_1}(t)d_\al(x^{m_1}y^{m_2}z^{m_3}))
\eqno(2.75)$$ is the solution space in $\mbb{C}[t,x,y,z]$ of the
equation (2.64) by Lemmas 2.1 and 2.4.\psp

{\bf Theorem 2.6}. {\it The followings are real solutions of the
equation (2.63): $$u=\cos
at\:\sum_{\al\in\mbb{N}^{11}}\zeta_{|\al|_1,0}d_\al(x^{m_1}y^{m_2}z^{m_3})-
\sin
at\:\sum_{\al\in\mbb{N}^{11}}\zeta_{|\al|_1,1}d_\al(x^{m_1}y^{m_2}z^{m_3}),
\eqno(2.76)$$
$$u=\sin
at\:\sum_{\al\in\mbb{N}^{11}}\zeta_{|\al|_1,0}d_\al(x^{m_1}y^{m_2}z^{m_3})+
\cos
at\:\sum_{\al\in\mbb{N}^{11}}\zeta_{|\al|_1,1}d_\al(x^{m_1}y^{m_2}
z^{m_3}), \eqno(2.77)$$ where $m_1,m_2,m_3\in\mbb{N}$.}\psp

The following lemma will be used to handle some special cases when
the operator $T_1$  in Lemma 2.1 does not have a right inverse and
to solve certain initial value problems in next section.\psp

 {\bf Lemma 2.7}. {\it Suppose that ${\cal A}$ is a free module
over a subalgebra ${\cal B}$ generated by a filtrated subspace
$V=\bigcup_{r=0}^\infty V_r$ (i.e., $V_r\subset V_{r+1}$). Let
$T_0$ be a linear operator on ${\cal A}$ with right inverse
$T_0^-$ such that
$$T_0({\cal B}),T_0^-({\cal B})\subset {\cal B},\;\;T_0(\eta_1\eta_2)=
T_0(\eta_1)\eta_2\qquad\for\;\;\eta_1 \in {\cal B},\;\eta_2\in
V,\eqno(2.78)$$ and let $T_1,...,T_m$ be commuting linear
operators on ${\cal A}$ such that $T_i(V)\subset V$,
$$T_0T_i=T_iT_0,\qquad T_i(f\zeta)=fT_i(\zeta)
\qquad\mbox{\it for}\;\;i\in\ol{1,m},\;f\in{\cal
B},\;\zeta\in{\cal A}.\eqno(2.79)$$ If $T_1^m(h)=0$ with
$h\in{\cal B}$ and $g\in V$, then
$$u=\sum_{i=0}^\infty(\sum_{s=1}^m(T_0^-)^sT_s)^i(hg)=
\sum_{i_1,...,i_m=0}^\infty {i_1+\cdots+i_m\choose
i_1,...,i_m}(T_0^-)^{\sum_{s=1}^msi_s}(h)(\prod_{r=1}^mT_r^{i_r})(g)\eqno(2.80)$$
is a solution of the equation:
$$(T_0^m-\sum_{r=1}^m T_0^{m-i}T_i)(u)=0.\eqno(2.81)$$
Suppose
$$T_i(V_r)\subset V_{r-1}\qquad\mbox{\it
for}\;\;i\in\ol{1,m},\;r\in\mbb{N},\eqno(2.82)$$ where
$V_{-1}=\{0\}$. Then any polynomial solution of (2.81) is a linear
combinations of the solutions of the form (2.80).}

{\it Proof}. Note that
$$T_0^{m-i}=T_0^m
(T_0^-)^i\qquad\for\;\;i\in\ol{1,m}\eqno(2.83)$$ and
\begin{eqnarray*} \hspace{2cm}& &\sum_{i_1+\cdots+i_m=i+1}{i+1\choose
i_1,...,i_m}y_1^{i_1}\cdots y_m^{i_m}=(y_1+\cdots+y_m)^{i+1}\\
&=&\sum_{r=1}^m\sum_{i_1+\cdots+i_m=i}{i\choose
i_1,...,i_m}y_ry_1^{i_1}\cdots
y_m^{i_m}.\hspace{5.3cm}(2.84)\end{eqnarray*}
 Thus
\begin{eqnarray*} & &(T_0^m-\sum_{p=1}^m T_0^{m-p}T_p)\left[\sum_{i_1,...,i_m=0}^\infty
{i_1+\cdots+i_m\choose
i_1,...,i_m}(T_0^-)^{\sum_{s=1}^msi_s}(h)(\prod_{r=1}^mT_r^{i_r})(g)\right]\\
&=&\sum_{i_1,...,i_m\in\mbb{N};\:i_1+\cdots+i_m>0}
{i_1+\cdots+i_m\choose
i_1,...,i_m}T_0^m(T_0^-)^{\sum_{s=1}^msi_s}(h)(\prod_{r=1}^mT_r^{i_r})(g)\\
& &- \sum_{i_1,...,i_m=0}^\infty\sum_{p=1}^m
{i_1+\cdots+i_m\choose
i_1,...,i_m}T_0^m(T_0^-)^{i_p+\sum_{s=1}^msi_s}(h)(T_p\prod_{r=1}^mT_r^{i_r})(g)=0.
\hspace{1.4cm}(2.85)\end{eqnarray*}

Suppose that (2.82) holds.  Let  $u\in {\cal B}V_k\setminus{\cal
B}V_{k-1}$ be a solution of (2.81). Take a basis
$\{\phi_i+V_{k-1}\mid i\in I\}$ of $V_k/V_{k-1}$. Write
$$u=\sum_{i\in I}h_i\phi_i+u',\qquad h_i\in{\cal B},\;u'\in{\cal
B}V_{k-1}.\eqno(2.86)$$ Since
$$T_r(\phi_i)\in V_{k-1}\qquad\for\;\;i\in
I,\;r\in\ol{1,m}\eqno(2.87)$$ by (2.82), we have
$$(T_0^m-\sum_{r=1}^m T_0^{m-i}T_i)(u)\equiv \sum_{i\in
I}T_0^m(h_i)\phi_i\equiv 0\;\;(\mbox{mod}\;{\cal
B}V_{k-1}).\eqno(2.88)$$ Hence
$$T_0^m(h_i)=0\qquad\for\;\;i\in I.\eqno(2.89)$$
Now
$$u-\sum_{j\in I}\sum_{i_1,...,i_m=0}^\infty
{i_1+\cdots+i_m\choose
i_1,...,i_m}(T_0^-)^{\sum_{s=1}^msi_s}(h_j)(\prod_{r=1}^mT_r^{i_r})(\phi_j)\in
{\cal B}V_{k-1}\eqno(2.90)$$ is a solution of (2.81). By induction
on $k$, $u$ is a linear combinations of the solutions of the form
(2.80).$\qquad\Box$\psp

We remark that the above lemma does not imply Lemma 2.1 because
$T_1$ and $T_2$ in Lemma 2.1 may not commute.

 Let $d_1$ be a differential operator on
$\mbb{F}[x_1,x_2,...,x_r]$ and let $d_2$ be a locally nilpotent
differential operator on $V=\mbb{F}[x_{r+1},...,x_n]$. Set
$$V_i=\{f\in V\mid d_2^i(f)\neq
0,\;d_2^{i+1}(f)=0\}\qquad\for\;\;i\in\mbb{N}.\eqno(2.91)$$ Take
 a subset $\{\psi_{i,j}\mid i\in\mbb{N}+1,\;j\in I_i\}$ of $V$ such
 that $\{\psi_{i,j}+V_{i-1}\mid j\in I_i\}$ forms a basis of
 $V_i/V_{i-1}$ for $i\in\mbb{N}+1$.
 Fix $h\in \mbb{F}[x_1,...,x_r]$.\psp

{\bf Lemma 2.8}. {\it Let $i$ be a positive integer. Suppose that
$$u=\sum_{j\in I_i}f_j\psi_{i,j}+u'\in
\mbb{F}[x_1,x_2,...,x_n]\eqno(2.92)$$ with
$f_j\in\mbb{F}[x_1,...,x_r]$ and $d_2^i(u')=0$ is a solution of
the equation:
$$(d_1-hd_2)(u)=0.\eqno(2.93)$$
Then the system
$$\xi_0=f_j,\;\;d_1(\xi_{s+1})=h\xi_s\qquad
\for\;\;s\in\ol{0,i-1}\eqno(2.94)$$ has a solution
$\xi_1,...,\xi_i\in \mbb{F}[x_1,...,x_r]$ for each $j\in I_i$.}

{\it Proof}. Denote $V_{-1}=\{0\}$. Observe that if $\{f_j+V_p\mid
p\in J\}$  is a linearly independent subset of $V_{p+1}/V_p$, then
$\{d_2^s(f_j)+V_{p-s}\mid p\in J\}$ is a linearly independent
subset of $V_{p-s+1}/V_{p-s}$ for $s\in\ol{1,p+1}$ by (2.91). By
induction, we take a subset$\{\phi_{i-s,j}\mid j\in J_{i-s}\}$ of
$V_{i-s}$ for each $s\in\ol{1,i}$ such that
$$\{d^s_2(\psi_{i,j_1})+V_{i-s-1},d_2^{s-p}(\phi_{i-p,j_2})+V_{i-s-1}\mid
p\in\ol{1,s},\;j_1\in I_i,\;j_2\in J_{i-p}\}\eqno(2.95)$$ forms a
basis of $V_{i-s}/V_{i-s-1}$ for $s\in\ol{1,i}$. Denote
$${\cal
B}=\sum_{s=1}^i\sum_{p=0}^{i-s}\sum_{j\in
J_{i-s}}\mbb{F}[x_1,...,x_r]d_2^p(\phi_{i-s,j}).\eqno(2.96)$$

Now we write
$$u=\sum_{j\in
I_i}[f_j\psi_{i,j}+\sum_{s=1}^if_{s,j}d^s(\psi_{i,j})]+v,\qquad
v\in{\cal B},\;f_{s,j}\in\mbb{F}[x_1,...,x_r].\eqno(2.97)$$ Then
(2.93) becomes \begin{eqnarray*}& &\sum_{j\in
I_i}[d_1(f_j)\psi_{i,j}+(d_1(f_{1,j})-hf_j)
d_2(\psi_{i,j})+\sum_{s=2}^i(d_1(f_{s,j})-hf_{s-1,j})d^s(\psi_{i,j})]
\\ & &+(d_1-hd_2)(v)=0.\hspace{10.7cm}(2.98)\end{eqnarray*} Since $(d_1-hd_2)(v)\in{\cal B}$,
we have:
$$d_1(f_j)=0,\;\;d_1(f_{1,j})=hf_j,\;\;d_1(f_{s,j})=hf_{s-1,j}\eqno(2.99)$$
for $j\in I_i$ and $s\in\ol{2,i}$. So (2.94) has a solution
$\xi_1,...,\xi_i\in \mbb{F}[x_1,...,x_r]$ for each $j\in
I_i.\hspace{0.7cm}\Box$\psp

Set
$${\cal S}_0=\{f\in \mbb{F}[x_1,...,x_r]\mid
d_1(f)=0\}\eqno(2.100)$$ and
$${\cal S}_i=\{f_0\in {\cal S}_0\mid d_1(f_s)=hf_{s-1}\;\mbox{for
some}\;\;f_1,...,f_i\in\mbb{F}[x_1,...,x_r]\}\eqno(2.101)$$ for
$i\in\mbb{N}+1$. For each $i\in\mbb{N}+1$ and $f\in {\cal S}_i$,
we fix $\{\sgm_1(f),...,\sgm_i(f)\}\subset\mbb{F}[x_1,...,x_r]$
such that
$$d_1(\sgm_1(f))=f,\;\;d_1(\sgm_s(f))=\sgm_{s-1}(f)\qquad\for\;\;s\in\ol{2,i}.
\eqno(2.102)$$ Denote $\sgm_0(f)=f$. By the proof of Lemma 2.1 and
the above Lemma, we obtain: \psp

{\bf Lemma 2.9}. {\it The set
$$\sum_{i=0}^\infty\sum_{j\in
I_i}\sum_{f\in{\cal
S}_i}\mbb{F}(\sum_{s=0}^i\sgm_s(f)d^s_2(\psi_{i,j}))\eqno(2.103)$$
is the solution space of the equation (2.93) in
$\mbb{F}[x_1,x_2,...,x_n]$.}\psp

Let $\es\in\{1,-1\}$ and let $\lmd$ be a nonzero real number. Next
we want to find all the polynomial solutions of the equation:
$$u_{tt}+\frac{\lmd}{t}u_t-\es(u_{x_1x_1}+u_{x_2x_2}+\cdots+u_{x_nx_n})=0,
\eqno(2.104)$$ which is the generalized anisymmetrical Laplace
equation (cf. [A]) if $\es=-1$. Rewrite the above equation as:
$$tu_{tt}+\lmd u_t-\es t(u_{x_1x_1}+u_{x_2x_2}+\cdots+u_{x_nx_n})=0.
\eqno(2.105)$$
 Set
$$d=t\ptl_t^2+\lmd\ptl_t.\eqno(2.106)$$
Denote
$${\cal S}=\{f\in\mbb{F}[t]\mid d(f)=0\}.\eqno(2.107)$$
Note that
$$d(t^m)=m(\lmd+m-1)t^{m-1}\qquad\for\;\;m\in\mbb{N}.\eqno(2.108)$$
So
$${\cal
S}=\left\{\begin{array}{ll}\mbb{F}&\mbox{if}\;\lmd\not\in-(\mbb{N}+1),\\
\mbb{F}+\mbb{F}t^{-\lmd+1}&\mbox{if}\;\lmd\in-(\mbb{N}+1).\end{array}\right.
\eqno(2.109)$$ In particular, $t^{-\lmd}\not\in d(\mbb{F}[t])$ and
so $d$ does not have a right inverse when $\lmd$ is a negative
integer.

 Set
$$\phi_0(t)=1,\;\;\phi_i(t)=\frac{t^{2i}}{i!2^i\prod_{r=0}^{i-1}(\lmd+2r+1)}
\eqno(2.110)$$ for $i\in\mbb{N}+1$ and when $\lmd\neq
-1,-3,...,-(2i-1)$. When $\lmd\in -(\mbb{N}+1)$, we set
$$\psi_0=t^{1-\lmd},\;\;\psi_i=\frac{t^{2i+1-\lmd}}{2^ii!\prod_{r=1}^i(2r+1-\lmd)}\qquad
\for\;\;i\in\mbb{N}+1.\eqno(2.111)$$ Define
$$V=\mbb{F}[x_1,x_2,...,x_n],\;\;\Dlt_n=\sum_{r=1}^n\ptl_{x_r}^2,\;\;\Dlt_{2,n}=
\sum_{s=2}^n\ptl_{x_s}^2 \eqno(2.112)$$ and
$$V_i=\{f\in V\mid
\Dlt_n^i(f)=0\}\qquad\for\;\;i\in\mbb{N}+1.\eqno(2.113)$$ Observe
\begin{eqnarray*}\hspace{1cm}&
&\sum_{j_1,...,j_i=0}^\infty(-1)^{j_1+\cdots+j_i}{j_1+\cdots+j_i\choose
j_1,...,j_i}\prod_{r=1}^i \left[{i\choose r}t^r\right]^{j_r}=
\sum_{p=0}^\infty\left(-\sum_{s=1}^i{i\choose s}t^i\right)^p\\ &=&
\frac{1}{(1+t)^i}=\sum_{r=0}^\infty (-1)^r{i+r-1\choose
r}t^r\hspace{6.6cm}(2.114)\end{eqnarray*} for $|t|<1.$
 Applying Lemma 2.7 to
$\Dlt_n^i=\sum_{r=0}^i{i\choose r}\ptl_{x_1}^{2(i-r)}\Dlt_{2,n}^r$
with $m=i$, $T_0=\ptl_{x_1}^2$ and $T_r=-{i\choose r}\Dlt_{2,n}^r$
for $r\in\ol{1,i}$, we get a basis
$$\left\{\sum_{r=0}^\infty (-1)^r{i+r-1\choose
r}\frac{x_1^{\ell_1+2r}}{(\ell_1+2r)!}\Dlt_{2,n}^r(x_2^{\ell_2}\cdots
x_n^{\ell_n})\mid
\ell_1\in\ol{0,2i-1},\;\ell_2,...,\ell_n\in\mbb{N}\right\}\eqno(2.115)$$
 of $V_i$. Hence we obtain:\psp

{\bf Theorem 2.10}. {\it If $\lmd\not\in -(\mbb{N}+1)$, then the
set
$$\{\sum_{r=0}^\infty
\es^r\phi_r(t)\Dlt_n^r(x_1^{\ell_1}\cdots x_n^{\ell_n})\mid
\ell_1,...,\ell_n\in\mbb{N}\}\eqno(2.116)$$ forms a basis of the
space of the polynomial solutions for the equation (2.104). When
$\lmd$ is a negative even integer, the set $$\{\sum_{r=0}^\infty
\es^r\phi_r(t)\Dlt_n^n(x_1^{\ell_1}\cdots
x_n^{\ell_n}),\sum_{r=0}^\infty
\es^r\psi_r(t)\Dlt_n^n(x_1^{\ell_1}\cdots x_n^{\ell_n})\mid
\ell_1,...,\ell_n\in\mbb{N}\}\eqno(2.117)$$ forms a basis of the
space of the polynomial solutions for the equation (2.104). Assume
that $\lmd=-2k-1$ is a negative odd integer. The set
\begin{eqnarray*}\hspace{1cm} & &\{\sum_{s=0}^k\sum_{r=0}^\infty
(-1)^r \es^s{k+r-1\choose r}\phi_s(t)\Dlt_n^s
\left[\frac{x_1^{\ell_1+2r}}{(\ell_1+2r)!}\Dlt_{2,n}^r(x_2^{\ell_2}\cdots
x_n^{\ell_n})\right],\\ & &\sum_{r=0}^\infty
\es^r\psi_r(t)\Dlt_n^r(x_1^{\ell_1'}x_2^{\ell_2}\cdots
x_n^{\ell_n}) \mid
\ell_1\in\ol{0,2i-1},\;\ell_1',\ell_2,...,\ell_n\in\mbb{N}\}\hspace{2cm}
(2.118)\end{eqnarray*} is a basis of the space of the polynomial
solutions for the equation (2.104).}\psp

Finally, we consider the special Euler-Poisson-Darboux equation:
$$u_{tt}-u_{x_1x_1}-u_{x_2x_2}-\cdots -u_{x_nx_n}
-\frac{m(m+1)}{t^2}u=0\eqno(2.119)$$ with $m\neq -1,0$. Change the
equations to:
$$t^2u_{tt}-t^2(u_{x_1x_1}+u_{x_2x_2}+\cdots +u_{x_nx_n})-m(m+1)u=0.\eqno(2.120)$$
Letting $u=t^{m+1}v$, we have:
$$t^2u_{tt}=m(m+1)t^{m+1}v+2(m+1)t^{m+2}v_t+t^{m+3}v_{tt}.\eqno(2.121)$$
Substituting (2.121) into (2.120), we get
$$tv_{tt}+2(m+1)v_t-t(v_{x_1x_1}+v_{x_2x_2}+\cdots +v_{x_nx_n})=0.
\eqno(2.122)$$ If we change variable $u=t^{-m}v$, then the
equation (2.120) becomes
$$tv_{tt}-2mv_t-t(v_{x_1x_1}+v_{x_2x_2}+\cdots+_{x_nx_n})=0
. \eqno(2.123)$$ Equations (2.121) and (2.122) are special cases
of the equation (2.105) with $\es=1$, and $\lmd=2(m+1)$ and
$\lmd=-2m$, respectively.

\section{Initial Value Problems}

In this section, we will solve two initial value problems.

Let $m$ and $n>1$ be positive integers and let
$$f_i(\ptl_{x_2},...,\ptl_{x_n})\in\mbb{R}[\ptl_{x_2},...,\ptl_{x_n}]
\qquad\for\;\;i\in\ol{1,m}.\eqno(3.1)$$ We want to solve the
equation:
$$(\ptl_{x_1}^m-\sum_{r=1}^m\ptl_{x_1}^{m-i}f_i(\ptl_{x_2},...,\ptl_{x_n}))(u)=0
\eqno(3.2)$$
with $x_1\in\mbb{R}$ and $x_r\in[-a_r,a_r]$ for $r\in\ol{2,n}$,
subject to the condition
$$\ptl_{x_1}^s(u)(0,x_2,...,x_n)=g_s(x_2,...,x_n)\qquad\for\;\;s\in\ol{0,m-1},
\eqno(3.3)$$
where $a_2,...,a_n$ are positive real numbers and
$g_0,...,g_{m-1}$ are continuous functions. For convenience, we
denote
$$k^\dg_i=\frac{k_i}{a_i},\;\;\vec
k^\dg=(k^\dg_2,...,k_n^\dg)\qquad\for\;\;\vec
k=(k_2,...,k_n)\in\mbb{N}^{\:n-1}.\eqno(3.4)$$ Set
$$e^{2\pi (\vec k^\dg\cdot\vec x)\sqrt{-1}}=e^{\sum_{r=2}^n2\pi
k^\dg_rx_r\sqrt{-1}}.\eqno(3.5)$$

For $r\in\ol{0,m-1}$, \begin{eqnarray*} & &
\frac{1}{r!}\sum_{i_1,...,i_m=0}^\infty {i_1+\cdots+i_m\choose
i_1,...,i_m}\int_{(x_1)}^{(\sum_{s=1}^msi_s)}(x_1^r)(\prod_{p=1}^mf_p(\ptl_{x_2},...,
\ptl_{x_n})^{i_p})(e^{2\pi (\vec k^\dg\cdot\vec x)\sqrt{-1}})\\
&=&\sum_{i_1,...,i_m=0}^\infty {i_1+\cdots+i_m\choose
i_1,...,i_m}\frac{x_1^{r+\sum_{s=1}^msi_s}}{(r+\sum_{s=1}^msi_s)!}
\\ & &\times\left[\prod_{p=1}^mf_p(2k_2^\dg\pi\sqrt{-1},...,2k_n^\dg\pi\sqrt{-1})^{i_p}\right]e^{2\pi
(\vec k^\dg\cdot\vec
x)\sqrt{-1}}\hspace{5.3cm}(3.6)\end{eqnarray*} is a complex
solution of the equation (3.2) by Lemma 2.7 for any $\vec
k\in\mbb{Z}^{\:n-1}$. We write
\begin{eqnarray*}& &\sum_{i_1,...,i_m=0}^\infty {i_1+\cdots+i_m\choose
i_1,...,i_m}\frac{x_1^r\prod_{p=1}^m(x_1^pf_p(2k_2^\dg\pi\sqrt{-1},...,2k_n^\dg\pi\sqrt{-1}))
^{i_p}} {(r+\sum_{s=1}^msi_s)!} \\&=&\phi_r(x_1,\vec
k)+\psi_r(x_1,\vec k)\sqrt{-1},\hspace{9.3cm}(3.7)\end{eqnarray*}
where $\phi_r(x_1,\vec k)$ and $\psi_r(x_1,\vec k)$ are real
functions. Moreover,
$$\ptl_{x_1}^s(\phi_r)(0,\vec k)=\dlt_{r,s},\;\;\ptl_{x_1}^s(\psi_r)(0,\vec
k)=0\qquad\for\;\;s\in\ol{0,r}.\eqno(3.8)$$ We define $\vec 0\prec
\vec k$ if its first nonzero coordinate is a positive integer. By
superposition principle and Fourier expansions, we get:\psp

{\bf Theorem 3.1}. {\it The solution of the equation (3.2) subject
to the condition (3.3) is:
\begin{eqnarray*}u&=&\sum_{r=0}^{m-1}\sum_{\vec 0\preceq\vec
k\in\mbb{Z}^{\:n-1}}[b_r(\vec k)(\phi_r(x_1,\vec k^\dg)\cos
2\pi(\vec k^\dg\cdot\vec x)-\psi_r(x_1,\vec k^\dg)\sin 2\pi(\vec
k^\dg\cdot\vec x))\\ & &+c_r(\vec k)(\phi_r(x_1,\vec k^\dg)\sin
2\pi(\vec k^\dg\cdot\vec x)+\psi_r(x_1,\vec k^\dg)\cos 2\pi(\vec
k^\dg\cdot\vec x))],\hspace{3.5cm}(3.9)\end{eqnarray*}with
\begin{eqnarray*}\hspace{1.6cm}b_r(\vec k)&=&\frac{1}{2^{n-2}a_2\cdots
a_n}\int_{-a_2}^{a_2}\cdots \int_{-a_n}^{a_n}g_r(x_2,...,x_n)\cos
2\pi (\vec k^\dg\cdot\vec x)\:dx_n\cdots dx_2\\&
&-\sum_{s=0}^{r-1}(b_s(\vec k)\ptl_{x_1}^r(\phi_s)(0,\vec
k)+c_s(\vec k)\ptl_{x_1}^r(\psi_s)(0,\vec k))
\hspace{3.5cm}(3.10)\end{eqnarray*} and
\begin{eqnarray*}\hspace{1cm}c_r(\vec k)&=&\frac{1}{2^{n-2}a_2\cdots
a_n}\int_{-a_2}^{a_2}\cdots \int_{-a_n}^{a_n}g_r(x_2,...,x_n)\sin
2\pi (\vec k^\dg\cdot\vec x)\:dx_n\cdots dx_2\\
&&-\sum_{s=0}^{r-1}(c_s(\vec k) \ptl_{x_1}^r(\phi_s)(0,\vec
k)-b_s(\vec k) \ptl_{x_1}^r(\psi_s)(0,\vec
k)).\hspace{3.9cm}(3.11)\end{eqnarray*} The convergence of the
series (3.9) is guaranteed  by the Kovalevskaya Theorem on the
existence and uniqueness of the solution of linear partial
differential equations when the functions in (3.3) are analytic.}

\pse

{\bf Remark 3.2}. (1) If we take $f_i=b_i$ with $i\in\ol{1,m}$ to
be constant functions and $\vec k=\vec 0$ in (3.6), we get $m$
fundamental solutions
$$\vf_r(x)=\sum_{i_1,...,i_m=0}^\infty {i_1+\cdots+i_m\choose
i_1,...,i_m}\frac{x^r\prod_{p=1}^m(b_px^p) ^{i_p}}
{(r+\sum_{s=1}^msi_s)!},
 \qquad r\in\ol{0,m-1},\eqno(3.12)$$ of the constant-coefficient ordinary
differential equation
$$y^{(m)}-b_1y^{(m-1)}-\cdots-b_{m-1}y'-b_m=0.\eqno(3.13)$$
Given the initial conditions:
$$y^{(r)}(0)=c_r\qquad\for\;\;r\in\ol{0,m-1},\eqno(3.14)$$
we define $a_0=c_0$ and
$$a_r=c_r-\sum_{s=0}^{r-1}\sum_{i_1,...,i_{r-s}\in\mbb{N};\:\sum_{p=1}^rpi_p=r-s}{r-s\choose
i_1,...,i_{r-s}}a_sb_1^{i_1}\cdots b_{r-s}^{i_{r-s}}\eqno(3.15)$$
by induction on $r\in\ol{1,m-1}$. Now the solution of (3.13)
subject to the condition (3.14) is exactly
$$y=\sum_{r=0}^{m-1}a_r\vf_r(x).\eqno(3.16)$$
From the above results, it seems that the following functions
$${\cal Y}_r(y_1,...,y_m)=\sum_{i_1,...,i_m=0}^\infty {i_1+\cdots+i_m\choose
i_1,...,i_m}\frac{y_1 ^{i_1}y_2 ^{i_2}\cdots y_m ^{i_m}}
{(r+\sum_{s=1}^msi_s)!}\qquad\for\;\;r\in\mbb{N} \eqno(3.17)$$ are
important natural functions. Indeed,
$${\cal Y}_1(x)=e^x,\;\;{\cal Y}_0(0,-x)=\cos x,\;\;{\cal
Y}_1(0,-x)=\frac{\sin x}{x},\eqno(3.18)$$
$$\vf_r(x)=x^r{\cal
Y}_r(b_1x,b_2x^2,...,b_mx^m)\eqno(3.19)$$ and
\begin{eqnarray*}& &\phi_r(x_1,\vec x)+\psi_r(x_1,\vec x)\sqrt{-1}\\ &=&x_1^r{\cal
Y}_r(x_1f_1(2k_2^\dg\pi\sqrt{-1},...,2k_n^\dg\pi\sqrt{-1})),...,
x_1^mf_m(2k_2^\dg\pi\sqrt{-1},...,2k_n^\dg\pi\sqrt{-1}))\hspace{1.1cm}(3.20)
\end{eqnarray*}
for  $r\in\ol{0,m}$.

(2) We can solve the  initial value problem (3.2) and (3.3) with
the constant-coefficient differential operators
$f_i(\ptl_2,...,\ptl_n)$ replaced by variable-coefficient
differential operators
$\phi_i(\ptl_2,...,\ptl_{n_1})\psi_i(x_{n_1+1},...,x_n)$ for some
$2<n_1<n$ , where $\phi_i(\ptl_2,...,\ptl_{n_1})$ are polynomials
in $\ptl_2,...,\ptl_{n_1}$ and $\psi(x_{n_1+1},...,x_n)$ are
polynomials in $x_{n_1+1},...,x_n$. \psp

Let ${\cal T}=({\cal N},{\cal E})$ be a tree with $n$ nodes. Now
we consider the following {\it generalized wave equation}
$$u_{tt}-d_{\cal T}(u)=0\eqno(3.21)$$
(cf. (1.6)) with $t\in\mbb{R}$ and $x_r\in[-a_r,a_r]$ for
$r\in\ol{1,n}$ subject to the condition
$$u(0,x_1,...,x_n)=g_0(x_1,...,x_n),\;u_t(0,x_1,...,x_n)=g_1(x_1,...,x_n).
\eqno(3.22)$$ Denote $$\vec
x=(x_1,...,x_n),\;\;k^\dg_i=\frac{k_i}{a_i},\;\;\vec
k^\dg=(k^\dg_1,...,k_n^\dg)\qquad\for\;\;\vec
k=(k_1,...,k_n)\in\mbb{N}^{\:n}.\eqno(3.23)$$ For $\es\in\{0,1\}$,
$$\sum_{i=0}^\infty \int_{(t)}^{(2i)}(t^\es)d_{\cal T}^i(e^{2\pi (\vec k^\dg\cdot\vec
x)\sqrt{-1}})\eqno(3.24)$$ are solutions of (3.21) by Lemma 2.1.
Moreover,
$$\sum_{i=0}^\infty \int_{(t)}^{(2i)}(1)d_{\cal T}^i=\frac{1}{2}(e^{td_{\cal
T}}+e^{-td_{\cal T}}),\;\;\sum_{i=0}^\infty
\int_{(t)}^{(2i)}(t)d_{\cal T}^i=\frac{1}{2}\int_{(t)}(e^{td_{\cal
T}}+e^{-td_{\cal T}}) .\eqno(3.25)$$

Recall that $\Psi$ is the set of all tips in ${\cal T}=({\cal
N},{\cal E})$. A node $\iota_j$ is called a {\it descendant} of
$\iota_i$ if $i<j$ and there exist a sequence
$$i_0=i<i_1<\cdots i_{r-1}<i_r=j\eqno(3.26)$$
such that
$$(\iota_{i_r},\iota_{r+1})\in{\cal
E}\qquad\for\;\;r\in\ol{0,r-1}.\eqno(3.27)$$ Set
$${\cal D}_i=\mbox{the set of all descendants
of}\;\iota_i.\eqno(3.28)$$ Let
$$\td\xi_r(t)=t\ptl_{x_r}^{m_r}\qquad\for\;\;\iota_r\in\Psi.\eqno(3.29)$$
Suppose that we have defined $\{\td\xi_s(t)\mid \iota_s\in{\cal
D}_i\}$.  Denote
$$\Theta_i=\{\iota_s\in{\cal N}\mid (\iota_i,\iota_s)\in{\cal
E}\}\subset{\cal D}_i.\eqno(3.30)$$ Now we define
$$\td\xi_i(t)=\int_0^t(\ptl_{x_i}+\sum_{\iota_s\in\Theta_i}\td\xi_s(y_i))^2
dy_i.\eqno(3.31)$$ By induction, we have defined all
$\{\td\xi_1(t),...,\td\xi_n(t)\}$. Moreover, we let
$$\xi_1(t,\ptl_{\xi_1},...,\ptl_{x_n})=\td\xi_1(t),\;\;\xi_i(t,\ptl_{\xi_1},...,\ptl_{x_n})
=x_{p(i)}\td\xi_i(t) \eqno(3.32)$$ for $i\in\ol{2,n}$, where
$\iota_{p(i)}$ is the unique (parent) node such that
$(\iota_{p(i)},\iota_i)\in {\cal E}$. According to (5.48) in [X2],
$$e^{td_{\cal T}}=e^{\xi_n(t,\ptl_{\xi_1},...,\ptl_{x_n})}e^{\xi_{n-1}(t,\ptl_{\xi_1},...,\ptl_{x_n})}
\cdots e^{\xi_1(t,\ptl_{\xi_1},...,\ptl_{x_n})}.\eqno(3.33)$$ In
fact,
$$e^{td_{\cal T}}(e^{2\pi (\vec k^\dg\cdot\vec x)\sqrt{-1}})
=e^{\sum_{i=1}^n\xi_i(t, 2\pi k^\dg_1\sqrt{-1},...,2\pi
k^\dg_n\sqrt{-1})}e^{2\pi (\vec k^\dg\cdot\vec
x)\sqrt{-1}}\eqno(3.34)$$  Define
\begin{eqnarray*}\phi_{\vec
k}(t,x_1,...,x_n) &=&\frac{1}{4}(e^{2\pi \vec (k^\dg\cdot\vec
x)\sqrt{-1}+\sum_{i=1}^n\xi_i(t, 2\pi k^\dg_1\sqrt{-1},...,2\pi
k^\dg_n\sqrt{-1})}\\ & &+e^{2\pi \vec (k^\dg\cdot\vec
x)\sqrt{-1}-\sum_{i=1}^n\xi_i(t, 2\pi k^\dg_1\sqrt{-1},...,2\pi
k^\dg_n\sqrt{-1})}\\ &&+e^{-2\pi (\vec k^\dg\cdot\vec
x)\sqrt{-1}+\sum_{i=1}^n\xi_i(t, -2\pi k^\dg_1\sqrt{-1},...,-2\pi
k^\dg_n\sqrt{-1})}\\ & &+e^{-2\pi \vec (k^\dg\cdot\vec
x)\sqrt{-1}-\sum_{i=1}^n\xi_i(t, -2\pi k^\dg_1\sqrt{-1},...,-2\pi
k^\dg_n\sqrt{-1})}),\hspace{3.3cm}(3.35)\end{eqnarray*}
\begin{eqnarray*}\psi_{\vec
k}(t,x_1,...,x_n)&=&\frac{1}{4\sqrt{-1}}(e^{2\pi \vec
(k^\dg\cdot\vec x)\sqrt{-1}+\sum_{i=1}^n\xi_i(t, 2\pi
k^\dg_1\sqrt{-1},...,2\pi k^\dg_n\sqrt{-1})}\\ & &+e^{2\pi \vec
(k^\dg\cdot\vec x)\sqrt{-1}-\sum_{i=1}^n\xi_i(t, 2\pi
k^\dg_1\sqrt{-1},...,2\pi k^\dg_n\sqrt{-1})}\\ &&-e^{-2\pi (\vec
k^\dg\cdot\vec x)\sqrt{-1}+\sum_{i=1}^n\xi_i(t, -2\pi
k^\dg_1\sqrt{-1},...,-2\pi k^\dg_n\sqrt{-1})}\\& &-e^{-2\pi \vec
(k^\dg\cdot\vec x)\sqrt{-1}-\sum_{i=1}^n\xi_i(t, -2\pi
k^\dg_1\sqrt{-1},...,-2\pi
k^\dg_n\sqrt{-1})})\hspace{3.4cm}(3.36)\end{eqnarray*} for $\vec
k\in\mbb{Z}^{\:n}$. Then
$$\phi_{\vec
k}(0,x_1,...,x_n)=\cos 2\pi (\vec k^\dg\cdot\vec x),\;\;\psi_{\vec
k}(0,x_1,...,x_n)=\sin 2\pi (\vec k^\dg\cdot\vec x),\eqno(3.37)$$
$$\ptl_t(\phi_{\vec
k})(0,x_1,...,x_n)=\ptl_t(\psi_{\vec
k})(0,x_1,...,x_n)=0.\eqno(3.38)$$ Again we define $\vec 0\prec
\vec k$ if its first nonzero coordinate is a positive integer. By
superposition principle and Fourier expansions, we obtain:\psp

{\bf Theorem 3.2}. {\it The solution of the equation (3.21)
subject to (3.22) is \begin{eqnarray*}\hspace{1cm}u&=&\sum_{\vec
0\preceq\vec k\in\mbb{Z}^{\:n}}[b_{0,\vec k}\phi_{\vec
k}(t,x_1,...,x_n)+c_{0,\vec k}\psi_{\vec k}(t,x_1,...,x_n)\\ & &+
b_{1,\vec k}\int_{(t)}\phi_{\vec k}(t,x_1,...,x_n)+c_{1,\vec
k}\int_{(t)}\psi_{\vec
k}(t,x_1,...,x_n)]\hspace{3.7cm}(3.39)\end{eqnarray*}
 with
$$b_{\es,\vec k}=\frac{1}{a_1a_2\cdots a_n}\int_{-a_1}^{a_1}\cdots
\int_{-a_n}^{a_n}g_\es(x_1,...,x_n)\cos 2\pi (\vec k^\dg\cdot\vec
x)\:dx_n\cdots dx_1\eqno(3.40)$$ and
$$c_{\es,\vec k}=\frac{1}{a_1a_2\cdots a_n}\int_{-a_1}^{a_1}\cdots
\int_{-a_n}^{a_n}g_\es(x_1,...,x_n)\sin 2\pi (\vec k^\dg\cdot\vec
x)\:dx_n\cdots dx_1.\eqno(3.41)$$ The convergence of the series
(3.39) is guaranteed  by the Kovalevskaya Theorem on the existence
and uniqueness of the solution of linear partial differential
equations when the functions in (3.22) are analytic.}

\pse

{\bf Example 3.1}. Consider the special case
$$u_{tt}-(\ptl_{x_1}^2+x_1\ptl_{x_2}^2+x_2\ptl_{x_3}^2)(u)=0\eqno(3.42)$$
of (3.21). By (2.71) and (2.72), we have
\begin{eqnarray*}& &\xi_1(t, 2\pi
k^\dg_1\sqrt{-1},2\pi k^\dg_12\sqrt{-1},2\pi
k^\dg_3\sqrt{-1})\\&=&
-4\pi^2t\left[(k_1^\dg)^2-\frac{4\pi^2t^2}{3}((k^\dg_2)^4+2k_1^\dg
k_2^\dg(k_3^\dg)^2)+\frac{16\pi^4t^4}{3}(k_2^\dg)^2(k_3^\dg)^4-\frac{64\pi^6t^6}{21}
(k_3^\dg)^8\right]\\ & &-8\pi^3
t^2\left[k_1^\dg(k_2^\dg)^2-\frac{2\pi^2t^2}{3}(3(k_2^\dg)^3(k_3^\dg)^2+k_1^\dg(k_3^\dg)^4)
+\frac{16t^4}{9}k_2^\dg(k_3^\dg)^6\right]\sqrt{-1},\hspace{1.6cm}(3.43)\end{eqnarray*}
\begin{eqnarray*}\hspace{2cm}& &\xi_2(t, 2\pi k^\dg_1\sqrt{-1},2\pi k^\dg_12\sqrt{-1},2\pi
k^\dg_3\sqrt{-1})\\
&=&-4\pi^2tx_1\left[(k_2^\dg)^2-\frac{4\pi^2t^2}{3}(k_3^\dg)^4\right]-8\pi^3k_2^\dg
(k_3^\dg)^3t^2x_1\sqrt{-1}\hspace{2.6cm}(3.44)\end{eqnarray*} and
$$\xi_3(t, 2\pi
k^\dg_1\sqrt{-1},2\pi k^\dg_12\sqrt{-1},2\pi k^\dg_3\sqrt{-1})
=-4\pi^2(k_3^\dg)^2tx_2.\eqno(3.45)$$ Thus
\begin{eqnarray*}& &\phi_{\vec
k}(t,x_1,x_2,x_3)=\\ & &\frac{1}{2}\exp
[4\pi^2t\left[(k_1^\dg)^2-\frac{4\pi^2t^2}{3}((k^\dg_2)^4+2k_1^\dg
k_2^\dg(k_3^\dg)^2)+\frac{16\pi^4t^4}{3}(k_2^\dg)^2(k_3^\dg)^4-\frac{64\pi^6t^6}{21}
(k_3^\dg)^8\right]\\ &
&+4\pi^2tx_1\left[(k_2^\dg)^2-\frac{4\pi^2t^2}{3}(k_3^\dg)^4\right]+
4\pi^2(k_3^\dg)^2tx_2]\cos [2\pi (\vec k^\dg\cdot\vec x)+8\pi^3
t^2[k_1^\dg(k_2^\dg)^2\\
& &-\frac{2\pi^2t^2}{3}(3(k_2^\dg)^3(k_3^\dg)^2
+k_1^\dg(k_3^\dg)^4)
+\frac{16t^4}{9}k_2^\dg(k_3^\dg)^6]+8\pi^3k_2^\dg
(k_3^\dg)^3t^2x_1]\\ & &+\frac{1}{2}\exp
[-4\pi^2t\left[(k_1^\dg)^2-\frac{4\pi^2t^2}{3}((k^\dg_2)^4+2k_1^\dg
k_2^\dg(k_3^\dg)^2)+\frac{16\pi^4t^4}{3}(k_2^\dg)^2(k_3^\dg)^4-\frac{64\pi^6t^6}{21}
(k_3^\dg)^8\right]\\
&
&-4\pi^2tx_1\left[(k_2^\dg)^2-\frac{4\pi^2t^2}{3}(k_3^\dg)^4\right]-
4\pi^2(k_3^\dg)^2tx_2]\cos [2\pi (\vec k^\dg\cdot\vec x)-8\pi^3
t^2[k_1^\dg(k_2^\dg)^2\\
& &-\frac{2\pi^2t^2}{3}(3(k_2^\dg)^3(k_3^\dg)^2
+k_1^\dg(k_3^\dg)^4)
+\frac{16t^4}{9}k_2^\dg(k_3^\dg)^6]-8\pi^3k_2^\dg
(k_3^\dg)^3t^2x_1]\hspace{3cm}(3.46)\end{eqnarray*}
\begin{eqnarray*}& &\psi_{\vec
k}(t,x_1,x_2,x_3)=\\ & &\frac{1}{2}\exp
[4\pi^2t\left[(k_1^\dg)^2-\frac{4\pi^2t^2}{3}((k^\dg_2)^4+2k_1^\dg
k_2^\dg(k_3^\dg)^2)+\frac{16\pi^4t^4}{3}(k_2^\dg)^2(k_3^\dg)^4-\frac{64\pi^6t^6}{21}
(k_3^\dg)^8\right]\\ &
&+4\pi^2tx_1\left[(k_2^\dg)^2-\frac{4\pi^2t^2}{3}(k_3^\dg)^4\right]+
4\pi^2(k_3^\dg)^2tx_2]\sin [2\pi (\vec k^\dg\cdot\vec x)+8\pi^3
t^2[k_1^\dg(k_2^\dg)^2\\
& &-\frac{2\pi^2t^2}{3}(3(k_2^\dg)^3(k_3^\dg)^2
+k_1^\dg(k_3^\dg)^4)
+\frac{16t^4}{9}k_2^\dg(k_3^\dg)^6]+8\pi^3k_2^\dg
(k_3^\dg)^3t^2x_1]\\ & &+\frac{1}{2}\exp
[-4\pi^2t\left[(k_1^\dg)^2-\frac{4\pi^2t^2}{3}((k^\dg_2)^4+2k_1^\dg
k_2^\dg(k_3^\dg)^2)+\frac{16\pi^4t^4}{3}(k_2^\dg)^2(k_3^\dg)^4-\frac{64\pi^6t^6}{21}
(k_3^\dg)^8\right]\hspace{3cm}\end{eqnarray*}\begin{eqnarray*} &
&-4\pi^2tx_1\left[(k_2^\dg)^2-\frac{4\pi^2t^2}{3}(k_3^\dg)^4\right]-
4\pi^2(k_3^\dg)^2tx_2]\sin [2\pi (\vec k^\dg\cdot\vec x)-8\pi^3
t^2[k_1^\dg(k_2^\dg)^2\\
& &-\frac{2\pi^2t^2}{3}(3(k_2^\dg)^3(k_3^\dg)^2
+k_1^\dg(k_3^\dg)^4)
+\frac{16t^4}{9}k_2^\dg(k_3^\dg)^6]-8\pi^3k_2^\dg
(k_3^\dg)^3t^2x_1]\hspace{3cm}(3.47)\end{eqnarray*}
 for
$\vec k\in\mbb{Z}^{\;3}$.

\section{Polynomial Representations of Lie Algebras}

In this section, we give three simple examples of applying Lemma
2.1 to polynomial representations of Lie algebras.

Let $\mbb{F}$ be any field with characteristic 0 and let $n\geq 3$
be an integer. The special orthogonal Lie algebra
$$so(n,\mbb{F})=\sum_{1\leq i<j\leq
n}\mbb{F}(E_{i,j}-E_{j,i}),\eqno(4.1)$$ where $E_{r,s}$ is an
$n\times n$ matrix with $1$ as its $(r,s)$-entry and 0 as the
others. There is a natural representation of $so(n,\mbb{F})$ on
the algebra ${\cal A}=\mbb{F}[x_1,...,x_n]$ of polynomials in $n$
variables:
$$(E_{i,j}-E_{j,i})|_{\cal A}=x_i\ptl_{x_j}-x_j\ptl_{x_i}\qquad\for\;\;1\leq i<j\leq
n.\eqno(4.2)$$ Denote
$$|\al|=\sum_{i=1}^n\al_i\qquad\for\;\;\al=(\al_1,...,\al_n)\in\mbb{N}^{\:n}\eqno(4.3)$$
and
$${\cal
A}_k=\sum_{\al\in\mbb{N}^{\:n},\;|\al|=k}\mbb{F}x^\al\qquad\for\;\;k\in\mbb{N},\eqno(4.4)$$
the space of homogeneous polynomials of degree $k$. Set
$${\cal H}_k=\{f\in {\cal A}_k\mid
(\ptl_{x_1}^2+\cdots+\ptl_{x_n}^2)(f)=0\},\eqno(4.5)$$ the space
of homogeneous harmonic polynomials of degree $k$. It is well
known in harmonic analysis that ${\cal H}_k$ are irreducible
$so(n,\mbb{F})$-submodules and
$${\cal A}_k={\cal H}_k\oplus (x_1^2+\cdots +x_n^2){\cal
A}_{k-2}\qquad\for\;\;2\leq k\in\mbb{N}.\eqno(4.6)$$ By (2.35), we
have: \psp

{\bf Theorem 4.1}. {\it The following set
\begin{eqnarray*}\hspace{1.7cm}& &\{\sum_{r_2,...,r_k=0}^{\infty}\frac{(-1)^{r_2+\cdots+
r_k}{r_2+\cdots+r_k\choose r_2,...,r_n}\prod_{i=2}^k{\ell_i\choose
2r_i}}{(1+2\es(r_2+\cdots+r_k)){2(r_2+\cdots+r_k)\choose
2r_2,...,2r_n}}x_1^{\es+2(r_2+\cdots +r_n)}\\
& &\times x_2^{\ell_2-2r_2}\cdots x_n^{\ell_n-2r_n}
\mid\es\in\{0,1\};\;\ell_2,...,\ell_n\in\mbb{N},\;\es+\sum_{i=2}\ell_i=k\}
\hspace{1.8cm}(4.5)\end{eqnarray*} forms a basis of ${\cal H}_k$.}
\psp

Recall the special linear Lie algebra
$$sl(n,\mbb{F})=\sum_{i\neq
j}\mbb{F}E_{i,j}+\sum_{i=1}^{n-1}\mbb{F}(E_{i,i}-E_{i+1,i+1}).\eqno(4.8)$$
Note
$$H=\sum_{i=1}^{n-1}\mbb{F}(E_{i,i}-E_{i+1,i+1})\eqno(4.9)$$
is a Cartan subalgebra of $sl(n,\mbb{F})$. Take $\{E_{i,j}\mid
1\leq i<j\leq n\}$ as positive root vectors. Let
$${\cal Q}=\mbb{F}(x_1,...,x_n,y_1,...,y_n),\eqno(4.10)$$
the space of rational functions in $x_1,...,x_n,y_1,...,y_n$.
Define a representation of $sl(n,\mbb{F})$ on ${\cal Q}$ via
$$E_{i,j}|_{\cal
Q}=x_i\ptl_{x_j}-y_j\ptl_{y_i}\qquad\for\;\;i,j\in\ol{1,n}.\eqno(4.11)$$
A nonzero function $f\in {\cal Q}$ is called  {\it singular} if
$$E_{i,j}(f)=0\qquad\for\;\;1\leq i<j\leq n\eqno(4.12)$$
and there exist a linear function $\lmd$ on $H$ such that
$$h(f)=\lmd(h)f\qquad\for\;\;h\in H.\eqno(4.13)$$

Set
$$\zeta=\sum_{i=1}^n x_iy_i.\eqno(4.14)$$
Then
$$\xi(\zeta)=0\qquad\for\;\;\xi\in sl(n,\mbb{F}).\eqno(4.15)$$
\pse

{\bf Lemma 4.2}. {\it Any singular function in ${\cal Q}$ is a
rational function in $x_1,y_n,\zeta$}.

{\it Proof}. Let $f\in{\cal Q}$ be a singular function. We can
write $$f=g(x_1,...,x_{n-1},\zeta,y_1,...,y_n)\eqno(4.16)$$ as a
rational functions in $x_1,...,x_{n-1},\zeta,y_1,...,y_n$. By
(4.11), (4.12) and (4.15), we have:
$$E_{i,n}(f)=-y_n\ptl_{y_i}(g)=0\qquad\for\;\;i\in\ol{1,n-1},\eqno(4.17)$$
equivalently,
$$\ptl_{y_i}(g)=0\qquad\for\;\;i\in\ol{1,n-1}.\eqno(4.19)$$
Thus (4.11) and (4.19) imply
$$E_{1,i}(g)=x_1\ptl_{x_i}(g)=0\qquad\for\;\;i\in\ol{2,n-1},\eqno(4.20)$$
that is,
$$\ptl_{x_i}(g)=0\qquad\for\;\;i\in\ol{2,n-1},\eqno(4.21)$$
Therefore, $g$ is independent of $x_2,...,x_{n-1}$ and
$y_1,...,y_{n-1}.\qquad\Box$ \psp

Set
$${\cal
A}_{\ell_1,\ell_2}=\sum_{\al,\be\in\mbb{N}^{\:n};\;|\al|=\ell_1,\;|\be|=\ell_2}
\mbb{F}x^\al
y^\be\qquad\for\;\;\ell_1,\ell_2\in\mbb{N}.\eqno(4.22)$$ Then
${\cal A}_{\ell_1,\ell_2}$ is a finite-dimensional
$sl(n,\mbb{F})$-submodule by (4.11). Recall that the fundamental
weights $\lmd_1,..,\lmd_{n-1}$ are linear functions on $H$ such
that
$$\lmd_i(E_{j,j}-E_{j+1,j+1})=\dlt_{i,j}\qquad\for\;\;i,j\in\ol{1,n-1}.\eqno(4.23)$$
The function $x_1^{\ell_1}y_n^{\ell_2}$ is a singular function of
weight $\ell_1\lmd_1+\ell_2\lmd_{n-1}$. According to the above
lemma, any singular polynomial in ${\cal A}_{\ell_1,\ell_2}$ must
be of the form $ax_1^{\ell_1-i}y_n^{\ell_2-i}\zeta^i$ for some
$0\neq a\in\mbb{F}$ and $i\in\mbb{N}$. Define
$$V_{\ell_1,\ell_2}=\mbox{the submodule
generated by}\;x_1^{\ell_1}y_n^{\ell_2}.\eqno(4.24)$$ According to
Weyl's Theorem of completely reducibility, ${\cal
A}_{\ell_1,\ell_2}$ is a direct sum of its irreducible submodules,
which are generated by its singular polynomials. So
$V_{\ell_1,\ell_2}$ is an irreducible highest weight module with
the highest weight $\ell_1\lmd_1+\ell_2\lmd_{n-1}$ and
$${\cal A}_{\ell_1,\ell_2}=V_{\ell_1,\ell_2}\oplus \zeta
{\cal A}_{\ell_1-1,\ell_2-1}\eqno(4.25)$$ as a direct sum of two
$sl(n,\mbb{F})$-submodules, where we treat $V_{i,j}=\{0\}$ if
$\{i,j\}\not\subset\mbb{N}$.

Denote
$$\Dlt=\sum_{i=1}^n\ptl_{x_i}\ptl_{y_i}.\eqno(4.26)$$
It can be verified that
$$\xi\Dlt=\Dlt\xi\qquad\for\;\;\xi\in sl(n,\mbb{F}),\eqno(4.27)$$
as operators on ${\cal Q}$. Set
$${\cal H}_{\ell_1,\ell_2}=\{f\in {\cal A}_{\ell_1,\ell_2}\mid
\Dlt(f)=0\}\eqno(4.28)$$
 Since $\Dlt(x_1^{\ell_1}y_n^{\ell_2})=0$,
we have $$V_{\ell_1,\ell_2}\subset {\cal
H}_{\ell_1,\ell_2}\eqno(4.29)$$ by (4.24) and (4.27). On the other
hand,
$${\cal A}_{\ell_1,\ell_2}=\bigoplus_{i=0}^\infty
\zeta^iV_{\ell_1-i,\ell_2-i}\eqno(4.30)$$ by (4.25) and induction.
Note
$$\Dlt\zeta=n+\zeta\Dlt+\sum_{i=1}^n(x_i\ptl_{x_i}+y_i\ptl_{y_i})\eqno(4.31)$$
as operators on ${\cal Q}$. Thus
$$\Dlt(\zeta^i
g)=\sum_{r=1}^i(n+\ell_1+\ell_2-2r)(\zeta^{i-1}g)=i(n+\ell_1+\ell_2-i-1)
\zeta^{i-1}g \eqno(4.32)$$ for $i\in\mbb{N}+1,\;g\in
V_{\ell_1-i,\ell_2-i}.$ Hence
$${\cal H}_{\ell_1,\ell_2}\bigcap \zeta {\cal
A}_{\ell_1-1,\ell_2-1}=\{0\}.\eqno(4.33)$$ Therefore,
$$V_{\ell_1,\ell_2}={\cal H}_{\ell_1,\ell_2}\eqno(4.34)$$
by (4.25) and (4.29). Now Lemma 2.1 gives:\psp

{\bf Theorem 4.3}. {\it The set
\begin{eqnarray*}&
&\{\sum_{i_2,...,i_n=0}^\infty\frac{(-1)^{i_2+\cdots+i_n}\prod_{r=2}^n{m_r\choose
i_r}{l_r\choose i_r}}{{m+i_2+\cdots+i_n\choose
m}{i_2+\cdots+i_n\choose
i_2,...,i_n}}x_1^{m+i_2+\cdots+i_n}y_1^{i_2+\cdots+i_n}\prod_{r=2}^nx_r^{m_r-i_r}
y_r^{l_r-i_r},\\
& &
\sum_{i_2,...,i_n=0}^\infty\frac{(-1)^{i_2+\cdots+i_n}\prod_{r=2}^n{m_r'\choose
i_r}{l_r\choose i_r}}{{m'+i_2+\cdots+i_n\choose
m'}{i_2+\cdots+i_n\choose
i_2,...,i_n}}x_1^{i_2+\cdots+i_n}y_1^{m'+i_2+\cdots+i_n}\prod_{r=2}^nx_r^{m_r'-i_r}
y_r^{l_r'-i_r}\\
& &\mid
m,m',m_r,n_r\in\mbb{N};\;m+\sum_{r=2}^nm_r=\sum_{r=2}^nm'_r=\ell_1,\;
\sum_{r=2}^nl_r=m'+\sum_{r=2}^nl'_r=\ell_2\}\hspace{0.9cm}(4.35)\end{eqnarray*}
for a basis of $V_{\ell_1,\ell_2}$}.\psp

In the Lie algebra $sl(7,\mbb{F})$, we set
$$h_1=-2E_{2,2}+E_{3,3}+E_{4,4}+2E_{5,5}-E_{6,6}-E_{7,7},\;\;
h_2=E_{2,2}-E_{3,3}-E_{5,5}+E_{6,6},\eqno(4.36)$$
$$E_1=\sqrt{2}(E_{1,2}-E_{5,1})-E_{3,7}+E_{4,6},\;\; E_2=E_{2,3}-E_{6.5},\;\qquad\eqno(4.37)$$
$$E_3=[E_1,E_2]=\sqrt{2}(E_{1,3}-E_{6,1})+E_{2,7}-E_{4,5},\eqno(4.38)$$
$$E_4=[E_1,E_3]/2=\sqrt{2}(E_{1,7}-E_{4,1})+E_{6,2}-E_{5,3},\eqno(4.39)$$
$$E_5=[E_1,E_4]/3=E_{4,2}-E_{5,7},\;\;E_6=[E_5,E_2]=E_{4,3}-E_{6,7},\eqno(4.40)$$
$$F_1=\sqrt{2}(E_{2,1}-E_{1,5})-E_{7,3}+E_{6,4},\qquad F_2=E_{3,2}-E_{5,6},\eqno(4.41)$$
$$F_3=\sqrt{2}(E_{3,1}-E_{1,6})+E_{7,2}-E_{5,4},\qquad F_5=E_{2,4}-E_{7,5},\eqno(4.42)$$
$$F_4=\sqrt{2}(E_{7,1}-E_{1,4})+E_{2,6}-E_{3,5},\qquad F_6=E_{3,4}-E_{7,6}.\eqno(4.43)$$
Then the exceptional Lie algebra of type $G_2$ is the Lie
subalgebra
$${\cal G}_2=\mbb{C}h_1+\mbb{C}h_2+\sum_{i=1}^6(\mbb{C}E_i+\mbb{C}F_i)\eqno(4.44)$$
of $sl(7,\mbb{F})$ (cf. [H]). Its Cartan subalgebra
$$H=\mbb{C}h_1+\mbb{C}h_2.\eqno(4.45)$$
 We choose $\{E_1,E_2,...,E_6\}$ as positive root vectors.

Let $Q$ be the space of rational functions in $\{x_i\mid
i\in\ol{1,7}\}$ and define a representation of ${\cal G}_2$ on $Q$
via
$$E_{i,j}|_Q=x_i\ptl_{x_j}\qquad\for\;\;i,j\in\ol{1,7}.\eqno(4.46)$$
A nonzero function $f\in Q$ is called {\it singular} with respect
to ${\cal G}_2$ if
$$ E_i(f)=0\qquad\for\;\;i\in\ol{1,6}\eqno(4.47)$$
and (4.13) holds. Define
$$\eta=x_1^2+2x_2x_5+2x_3x_6+2x_4x_7.\eqno(4.48)$$
It can be verified that
$$\xi(\eta)=0\qquad\for\;\;\xi\in{\cal G}_2.\eqno(4.49)$$
\psp

{\bf Lemma 4.4}. {\it Any singular function in $Q$ with respect to
${\cal G}_2$ must be a rational function in $x_4$ and $\eta$.}

{\it Proof}. Let $f$ a singular function in $Q$. We can write
$$f=\vf(x_1,....,x_6,\eta)\eqno(4.50)$$
as a rational function in $x_1,...,x_6,\eta$. Note
$$E_5(f)=(x_4\ptl_{x_2}-x_5\ptl_{x_7})(\vf(x_1,....,x_6,\eta))
=x_4\ptl_{x_2}(\vf)=0\eqno(4.51)$$ and
$$E_6(f)=(x_4\ptl_{x_3}-x_6\ptl_{x_7})(\vf(x_1,....,x_6,\eta))
=x_4\ptl_{x_3}(\vf)=0.\eqno(4.52)$$ So
$$\ptl_{x_2}(\vf)=\ptl_{x_3}(\vf)=0.\eqno(4.53)$$

Next
$$E_2(f)=(x_2\ptl_{x_3}-x_6\ptl_{x_5})(\vf(x_1,....,x_6,\eta))
=-x_6\ptl_{x_5}(\vf)=0,\eqno(4.54)$$ that is,
$$\ptl_{x_5}(\vf)=0.\eqno(4.55)$$
Moreover,
$$E_3(f)=(\sqrt{2}(x_1\ptl_{x_3}-x_6\ptl_{x_1})
+x_2\ptl_{x_7}-x_4\ptl_{x_5})(\vf(x_1,....,x_6,\eta))
=-\sqrt{2}x_6\ptl_{x_1}(\vf)=0,\eqno(4.56)$$ which implies
$$\ptl_{x_1}(\vf)=0.\eqno(4.57)$$
Furthermore,
$$E_2(f)=(\sqrt{2}(x_1\ptl_{x_2}-x_5\ptl_{x_1})
-x_3\ptl_{x_7}+x_4\ptl_{x_6})(\vf(x_1,....,x_6,\eta))
=x_4\ptl_{x_6}(\vf)=0\eqno(4.58)$$ yields
$$\ptl_{x_6}(\vf)=0.\eqno(4.59)$$
According to (4.53), (4.55), (4.57) and (4.59), $\vf$ is
independent of $x_1,x_2,x_3,x_5,x_6.\qquad\Box$ \psp

Set
$$\Dlt'=\ptl_{x_2}^2+2\ptl_{x_2}\ptl_{x_5}+2\ptl_{x_3}\ptl_{x_6}+2\ptl_{x_4}\ptl_{x_7}.
\eqno(4.60)$$ It can be verified that
$$\Dlt'\xi=\xi\Dlt'\qquad\for\;\;\xi\in{\cal G}_2\eqno(4.60)$$
and
$$\Dlt'\eta=\eta\Dlt'+14+4\sum_{i=1}^7x_i\ptl_{x_i}\eqno(4.61)$$
as operators on $Q$. The fundamental weight $\lmd_1$ is a linear
function on $H$ such that
$$\lmd(h_1)=1,\qquad\lmd_1(h_2)=0.\eqno(4.62)$$
For $k\in\mbb{N}$, $x_4^k$ is a singular function with weight
$k\lmd_1$. Note that
$${\cal
B}_k=\sum_{\al\in\mbb{N}^{\:7};\;|\al|=k}\mbb{F}x^\al\eqno(4.63)$$
is a finite-dimensional ${\cal G}_2$-submodule. Let
$$V_k=\mbox{the submodule generated by}\;x_4^k.\eqno(4.64)$$
Then $V_k$ is an irreducible highest weight submodule of ${\cal
B}_k$ with highest weight $k\lmd_1$. By similar arguments as those
in the above of Theorem 4.3, we can prove
$$V_k=\{f\in{\cal B}_k\mid \Dlt'(f)=0\}.\eqno(4.65)$$
\pse

{\bf Theorem 4.5}. {\it The set
\begin{eqnarray*}\hspace{1cm}& &\{\sum_{i_2,i_3,i_4=0}^\infty\frac{(-1)^
{i_2+i_3+i_4}\prod_{r=2}^3{m_r\choose i_r}{m_{4+r}\choose
i_r}}{(1+2\es(i_2+i_3+i_4)){2(i_2+i_3+i_4)\choose
i_2,i_3,i_4,i_2+i_3+i_4}}x_1^{\es+2(i_2+i_3+i_4)}\prod_{s=2}^4x_s^{m_s-i_s}
x_{4+s}^{m_{4+s}-i_s}\\
&
&|\es\in\{0,1\},\;m_2,...,m_7\in\mbb{N};\;\es+\sum_{r=2}^7m_r=k\}
\hspace{5.2cm}(4.66)\end{eqnarray*} forms a basis of $V_k$}.

\vspace{1cm}

\noindent{\Large \bf References}

\hspace{0.5cm}

\begin{description}

\item[{[A]}]  A. V. Aksenov, Symmetries and fundamental solutions
of multidimensional generalized axi-symmetric equations Laplace
Equation, {\it Differentsial'nye Uravneniya} {\bf 29} (1993), 11-

\item[{[BG1]}] J. Barros-Neto and I. M. Gel'fand, Fundamental
solutions for the Tricomi operator, {\it Duke Math. J.} {\bf 98}
(1999), 465-483.

\item[{[BG2]}] J. Barros-Neto and I. M. Gel'fand, Fundamental
solutions for the Tricomi operator II, {\it Duke Math. J.} {\bf
111} (2002), 561-584.

\item[{[B]}] Yu. Yu. Berest, Weak invariants of local
transformation groups, {\it Differentsial'nye Uravneniya} {\bf 29}
(1993), 1796.

\item[{[I1]}] N. H. Ibragimov, {\it Transformation groups applied to
mathematical physics}, Nauka, 1983.

\item[{[I2]}] N. H. Ibragimov, {\it Lie Group Analysis of
Differential Equations}, Volume 2, CRC Handbook, CRC Press, 1995.

\item[{[H]}] J. E. Humphreys, {\it Introduction to Lie Algebras and Representation Theory},
 Springer-Verlag New York Inc., 1972.

\item[{[X1]}] X. Xu, Differential invariants of classical groups, {\it Duke Math. J.} {\bf 94}
(1998), 543-572.

\item[{[X2]}] X. Xu, Tree diagram Lie algebras of  differential
operators and evolution partial differential equations, {\it J.
Lie Theory} {\bf 16} (2006), no. 4, 691-718.

\end{description}

\end{document}